\documentclass[12pt,fleqn, a4]{article}

\usepackage[french]{babel}
\usepackage{inputenc}
\usepackage[T1]{fontenc}

\usepackage[mathscr]{euscript}
\usepackage{amsmath}

\usepackage{amssymb}
\usepackage{latexsym}
\usepackage{graphics}

\usepackage{graphics}
\usepackage{color}
\newcommand{\colval}{0.3}
\definecolor{colone}{gray}{\colval}

\newcommand{\dcb}{\begin{array}{lll}}
\newcommand{\dce}{\end{array}}
\newcommand{\ebe}{\begin{enumerate}\setlength{\baselineskip}{13pt}\setlength{\parskip}{5pt}}
\newcommand{\dbe}{\end{enumerate}}

\newcommand{\ibegin}{\begin{itemize}\setlength{\baselineskip}{19pt}\setlength{\parskip}{7pt}}
\newcommand{\iend}{\end{itemize}}
\newcommand{\ok}{\rule{4pt}{6pt}}
\newcommand{\desb}{\begin{description}}
\newcommand{\dese}{\end{description}}

\newtheorem{Thm}{Theorem}[section]
\newtheorem {Cor}[Thm]{Corollary}
\newtheorem {definition}[Thm]{Definition}
\newtheorem {pro}{Proposition}[Thm]
\newtheorem {Lemma}[Thm]{Lemma}
\newtheorem {rem}[Thm]{Remark}

\newtheorem {assumption}[Thm]{Assumption}

\newcommand {\bd}{\begin{definition}}
\newcommand {\ed}{\end{definition}}
\newcommand {\bpro}{\begin{pro}}
\newcommand {\epro}{\end{pro}}
\newcommand {\bl}{\begin{Lemma}}
\newcommand {\el}{\end{Lemma}}
\newcommand {\bcor}{\begin{Cor}}
\newcommand {\ecor}{\end{Cor}}
\newcommand {\brem }{\begin{rem} \rm }
\newcommand {\erem }{\end{rem}}
\newcommand{\bethe}{\begin{Thm}}
\newcommand{\ethe}{\end{Thm}}
\newcommand {\bassumption}{\begin{assumption}}
\newcommand {\eassumption}{\end{assumption}}

\newcommand{\transp}{{^\top\!}}

\def \ind{1\!\!1}
\def\cro#1{\langle #1\rangle}

\begin{document}

\begin{center}
\Large
Martingale representation processes and applications in the market viability with information flow expansion\footnote{A third version}
\end{center}

\begin{center}
Shiqi Song

{\footnotesize Laboratoire Analyse et Probabilités\\
Université d'Evry Val D'Essonne, France\\
shiqi.song@univ-evry.fr}
\end{center}

\

\begin{center}
Abstract
\end{center}

\begin{quote}
When the \textit{martingale representation property} holds, we call any local martingale which realizes the representation a \textit{representation process}. There are two properties of the \textit{representation process} which can greatly facilitate the computations under the \textit{martingale representation property}. Actually, on the one hand, the \textit{representation process} is not unique and there always exists a \textit{representation process} which is locally bounded and has pathwisely orthogonal components outside of a predictable thin set. On the other hand, the jump measure of a \textit{representation process} satisfies the \textit{finite predictable constraint}. In this paper, we give a detailed account of these two properties. As application, we will prove that, under the \textit{martingale representation property}, the \textit{full viability} of an expansion of market information flow implies the \textit{drift multiplier assumption}.
\end{quote}

\vspace{2cm}

\textbf{Key words.} Martingale representation property, enlargement of filtrations, hypothesis$(H')$, drift operator, market viability, local martingale deflator, \texttt{NA}$_{1}$ condition, conditional multiplicity, finite predictable constraint, martingale projection property.

\textbf{MSC class.} 60G07, 60G44, 60G40.

\

\section{Introduction}

The present paper, jointly with \cite{song-drift, songsuccessive}, takes part in a research program about the \textit{full viability} problem (cf. Definition \ref{fullvaibility}). Precisely, we consider a probability space $(\Omega,\mathcal{A},\mathbb{P})$ (with $\mathcal{A}$ a $\sigma$-algebra and $\mathbb{P}$ a probability measure on $\mathcal{A}$) and any pair consisting of a filtration $\mathbb{F}=(\mathcal{F}_{t})_{t\geq 0}$ of sub-$\sigma$-algebras of $\mathcal{A}$ and an $\mathbb{F}$ semimartingale $S$, which satisfies the \textit{no-arbitrage condition of the first kind} (\texttt{NA}$_{1}$ in abbreviation, cf. \cite{KC2010} for definition) in $\mathbb{F}$. Such a pair $(\mathbb{F},S)$ constitutes a model of financial market, where $\mathbb{F}$ represents the information flow and $S$ represents the asset process. Because of the no-arbitrage condition \texttt{NA}$_{1}$, the utility optimization problems in the model $(\mathbb{F},S)$ have solutions (cf. \cite{KC2010}) so that the model is also called a \textit{viable} market model. We are concerned with the consequences of a change in the information flow $\mathbb{F}$. (A better informed agent would operate with a market model with a bigger information flow $\mathbb{G}=(\mathcal{G}_{t})_{t\geq 0}$ ($\mathcal{F}_{t}\subset \mathcal{G}_{t}$) and he would manage his portfolio differently in regards to a less informed agent. How to quantify this difference would be essential for many purposes.) We are especially interested in those changes which preserve fully the market \textit{viability}, i.e., those filtrations $\mathbb{G}$$(\supset \mathbb{F})$ such that all the asset processes $S$ satisfying the \texttt{NA}$_{1}$ in $\mathbb{F}$ continue to satisfy \texttt{NA}$_{1}$ in $\mathbb{G}$ (the \textit{full viability} of $\mathbb{F}\subset\mathbb{G}$). A satisfactory result on this subject has been obtained in \cite{song-drift}.

A long computation has been necessary to establish the result in \cite{song-drift}, which has been carried out under various conditions, notably the \textit{martingale representation property} in $\mathbb{F}$ and the \textit{drift multiplier assumption} (cf. Definition \ref{assump1}). The former is well known, but the computation in \cite{song-drift} depends on features of the \textit{martingale representation property}, which are not so widely considered in the literature, namely the \textit{finite predictable constraint} (cf. Definition \ref{FPCC}) and the reconstruction of the \textit{representation processes}  (cf. Section \ref{mrp} for definition). As for the latter, the \textit{drift multiplier assumption} is a new notion in the literature of the theory of filtration enlargement. It has encountered serious questions about its relevance and its usefulness. As part of the research program, the present paper is specially devoted to the above mentioned properties.

The \textit{finite predictable constraint} is involved in our work because of the \textit{martingale projection property}. This latter property is defined for any (multi-dimensional) local martingale $M$ by the fact that, for any real local martingale $Y$ such that the (vector of) predictable dual projection $[Y,M]^{\mathbb{F}\cdot p}$ exists, there exists a predictable process $H$, integrable with respect to $M$, such that the predictable dual projection $[H\centerdot M,M]^{\mathbb{F}\cdot p}$ exists and the identity $$
[Y,M]^{\mathbb{F}\cdot p}= [\transp H\centerdot M,M]^{\mathbb{F}\cdot p}
$$ 
holds. In the sense of Lemma \ref{projectionlemma}, the \textit{martingale projection property} is equivalent to another formula:
$$
\{W{_*}(\mu-\nu): \mbox{ $W$ is $(\mu-\nu)$-integrable in the sense of \cite{Jacod}}\}
=
\{\transp H\centerdot M: \mbox{ $H$ is $M$-integrable}\},
$$
where $\mu$ denotes the jump measure of $M$ and $\nu$ denotes the $\mathbb{F}$ compensator of $\mu$. The \textit{martingale projection property} is needed, for example, in the proof of \cite[Theorem 4.3]{song-drift} or in the proof of Theorem \ref{driftsatisfied} below. The basic question, then, is how to recognize a local martingale $M$ which satisfies the \textit{martingale projection property}. It seems not very widely known, but the answer exists since Jacod \cite{Jacodlivre}, which consists to verify if the jump measure of $M$ satisfies the condition of \textit{finite predictable constraint}. (This last notion is therefore extracted from \cite[Th\'eor\`eme 4.80]{Jacodlivre} and a name is given to it because of its importance in the study of the \textit{full viability}.) For the applications in this paper and in \cite{song-drift}, we will extend the \textit{finite predictable constraint} condition to any integer valued random measure $\mu$ and we will give in Section \ref{sectionFPCC} a detailed account of the space$$
\{W{_*}(\mu-\nu): \mbox{ $W$ is $(\mu-\nu)$-integrable in the sense of \cite{Jacod}}\}
$$
under that condition. The main result is Theorem \ref{specif}. Specifications in the cases of accessible and respectively inaccessible \textit{time supports} (cf. Section \ref{irm} for definition) are given in Section \ref{pjpacc} and Section \ref{pjpinacc}.  

Usually the \textit{martingale representation property} is mentioned to characterize a specific process (a Brownian motion, for example). But, in this paper, what is relevant is a stochastic basis having a \textit{martingale representation property}, whatever the \textit{representation processes} are. Yet more, we should make use of various different \textit{representation processes} to make easier the computations under the \textit{martingale representation property}. In Section \ref{mrp}, based on the fact that the \textit{martingale representation property} in the filtration $\mathbb{F}$ implies the finite \textit{conditional multiplicity} condition of the filtration $\mathbb{F}$ (cf. \cite{Jacodlivre} and \cite{BEKSY} and Section \ref{multiplicity} Lemma \ref{partition}), we prove in Theorem \ref{pathortho} and Corollary \ref{XXXo} and Theorem \ref{boundedW} that, when the \textit{martingale representation property} holds, it is always possible to reconstruct the \textit{representation process} so that it becomes locally bounded and has pathwisely orthogonal components outside of a predictable thin set. We recall that a first reason to consider the notion of \textit{conditional multiplicity} comes from \cite[Th\'eor\`eme 4.80]{Jacodlivre}. Its name is borrowed from \cite{BEKSY}. An interesting application of this notion can be found in \cite{BEKSY} for a study of Brownian filtrations. 

By the way, we recall in Section \ref{rpfpcc} how the \textit{martingale representation property} in the filtration $\mathbb{F}$ also implies the \textit{finite predictable constraint} condition of the jump measure of the \textit{representation processes}. (It is this implication which makes possible the application of Theorem \ref{specif} in the computation of \cite{song-drift}.) 

In application of the reconstruction of \textit{representation processes}, we study in Section \ref{dma} the \textit{drift multiplier assumption}. We will prove in Theorem \ref{driftsatisfied} that, under the \textit{martingale representation property} in $\mathbb{F}$, the \textit{full viability} of an expansion of information flow $\mathbb{F}\subset \mathbb{G}$ implies the \textit{Hypothesis}$(H')$, and the corresponding \textit{drift operator} satisfies the \textit{drift multiplier assumption} (cf. Section \ref{dma} for definition). This means that every $\mathbb{F}$ local martingale $X$ is a $\mathbb{G}$ special semimartingale and there exist a common (multi-dimensional) $\mathbb{F}$ local martingale $N$ and a $\mathbb{G}$ predictable process $\varphi$ such that the drift part $\Gamma(X)$ of $X$ in $\mathbb{G}$ writes in the form
$$
{(\mbox{\scriptsize $\triangledown$})}\ \ \ \ \ \ \Gamma(X)=\transp\varphi \centerdot[N,X]^{\mathbb{F}\cdot p}.
$$
The main points of this result are (\textbf{a}) the $\mathbb{F}$ local martingale $N$ and the $\mathbb{G}$ predictable process $\varphi$ are common for all $\mathbb{F}$ local martingale $X$ and (\textbf{b}) the predictable dual projection $[N,X]^{\mathbb{F}\cdot p}$ exists for all $\mathbb{F}$ local martingale $X$. With a \textit{martingale representation property} in $\mathbb{F}$, this result can be proved in the following way. First of all, in this case, instead of considering a general $\mathbb{F}$ local martingale $X$ in formula (\mbox{\scriptsize $\triangledown$}), we only need to consider the components of a \textit{representation process} $(X_{k})$. Also, to define $N$, we only need to find its coefficients $H_{k}$ in its martingale representation $N=\sum_{k}H_{k}{{\centerdot}X_{k}}$. Consequently, formula (\mbox{\scriptsize $\triangledown$}) changes into an equation in the unknowns $\varphi_{k}$ and $H_{k}$:
$$
{(\mbox{\scriptsize $\triangledown\triangledown$})}\ \ \ \ \ \  \Gamma(X_{h})=\sum_{k} \transp\varphi_{k} H_{k}\centerdot[X_{k},X_{h}]^{\mathbb{F}\cdot p},
$$
for any components $X_{h}$ of the \textit{representation process}. The solution of equation (\mbox{\scriptsize $\triangledown\triangledown$}) (if exists) is not unique. We can restrict our search among the coefficients $H_{k}$ which are locally bounded. When, also, a locally bounded \textit{representation process} is chosen, the local martingale $N$ thus defined will be locally bounded so that the predictable dual projection $[N,X]^{\mathbb{F}\cdot p}$ exists for any local martingale $X$. Hence, to prove Theorem \ref{driftsatisfied}, it is enough to solve equation (\mbox{\scriptsize $\triangledown\triangledown$}). For that, we firstly modify the \textit{representation process} to render the equation as simple as possible, and we apply Lemma \ref{fbd} to link equation (\mbox{\scriptsize $\triangledown\triangledown$}) to the assumption of the \textit{full viability}. After these treatments, equation (\mbox{\scriptsize $\triangledown\triangledown$}) can be solved by elementary computations.

Section \ref{gcfv} presents some general consequences of the \textit{full viability} assumption on an expansion of information flow $\mathbb{F}\subset \mathbb{G}$. For example, we will see that the $\mathbb{F}$ inaccessible stopping times remain inaccessible in $\mathbb{G}$, or that $A^{\mathbb{G}\cdot p}$ is absolutely continuous with respect to $A^{\mathbb{F}\cdot p}$ for any $\mathbb{F}$ adapted locally bounded increasing process $A$. These results are interesting in themselves. We also emphasize the benefits of working with the \textit{drift multiplier assumption} in an enlarged filtration.

Notice that in this introduction we have written the key notions in italic. This rule will be left out in the rest of the paper. Notice also that the assumption of the \textit{martingale representation property} in $\mathbb{F}$ needs not mean that the market model $(\mathbb{F},S)$ is complete.

\

\section{Notation and convention}

We work on a probability space $(\Omega,\mathcal{A},\mathbb{P})$ endowed with a filtration $\mathbb{F}=(\mathcal{F}_{t})_{t\geq 0}$ of sub-$\sigma$-algebras of $\mathcal{A}$, satisfying the usual conditions. We employ the vocabulary of stochastic calculus as defined in \cite{HWY, Jacodlivre} with the specifications below.

Relations between random variables is to be understood almost sure relations. For a random variable $X$ and a $\sigma$-algebra $\mathcal{F}$, the expression $X\in\mathcal{F}$ means that $X$ is $\mathcal{F}$-measurable. The notation $\mathbf{L}^p(\mathbb{P},\mathcal{F})$ denotes the space of $p$-times $\mathbb{P}$-integrable $\mathcal{F}$-measurable random variables.

By definition, $\Delta_0X=0$ for any c\` adl\`ag process $X$. A process $A$ with finite variation considered in this paper is automatically assumed c\`adl\`ag. We denote by $\mathsf{d}A$ the (signed) random measure that $A$ generates.

Different vector spaces $\mathbb{R}^d$ are used in the paper. An element $v$ in $\mathbb{R}^d$ is considered as a vertical vector. We denote its transposition by $\transp v$. We denote (indifferently) the null vector by $\boldsymbol{0}$.

We deal with finite family of real processes $X=(X_i)_{1\leq i\leq d}$ ($d\in\mathbb{N}^*$). It will be considered as process taking values in the vector space $\mathbb{R}^d$. To mention such an $X$, we say that $X$ is a $d$-dimensional process. In general we denote by $X_{i}$ the $i$th component of the vector $X$. When $X$ is a semimartingale, we denote by $[X,\transp X]$ the $d\times d$-dimensional matrix valued process whose components are $[X_i,X_j]$ for $1\leq i,j\leq k$.

With respect to the filtration $\mathbb{F}$, the notation ${^{\mathbb{F}\cdot p}}\bullet$ denotes the predictable projection, and the notation $\bullet^{\mathbb{F}\cdot p}$ denotes the predictable dual projection.

For any $\mathbb{F}$ special semimartingale $X$, we can decompose $X$ in the form (see \cite[Theorem 7.25]{HWY}) :$$
\dcb
X=X_0+X^m+X^v,\
X^m=X^c+X^{da}+X^{di},
\dce
$$
where $X^m$ is the martingale part of $X$ and $X^v$ is the predictable part of finite variation of $X$, $X^c$ is the continuous martingale part, $X^{da}$ is the part of compensated sum of accessible jumps, $X^{di}$ is the part of compensated sum of totally inaccessible jumps. We recall that this decomposition of $X$ depends on the reference probability and the reference filtration. In the computations below we apply this notation system only for the decompositions in $\mathbb{F}$. We recall that every part of the decomposition of $X$, except $X_0$, is assumed null at $t=0$.

In this paper we employ the notion of stochastic integral only about the predictable processes. The stochastic integral are defined as 0 at $t=0$. We use a point "$\centerdot$" to indicate the integrator process in a stochastic integral. For example, the stochastic integral of a real predictable process ${H}$ with respect to a real semimartingale $Y$ is denoted by ${H}\centerdot Y$, while the expression $\transp{K}(\centerdot[X,\transp X]){K}$ denotes the process$$
\int_0^t \sum_{i=1}^k\sum_{j=1}^k({K}_s)_{i,s}({K}_s)_{j,s} \mathsf{d}[X_i,X_j]_s,\ t\geq 0,
$$
where ${K}$ is a $k$-dimensional predictable process and $X$ is a $k$-dimensional semimartingale. The expression $\transp{K}(\centerdot[X,\transp X]){K}$ respects the matrix product rule. The value at $t\geq 0$ of a stochastic integral will be denoted, for example, by $\transp{K}(\centerdot[X,\transp X]){K}_t$.

The notion of the stochastic integral with respect to a $d$-dimensional local martingale $X$ follows \cite{Jacodlivre}. We say that a $d$-dimensional $\mathbb{F}$ predictable process is integrable with respect to $X$ under the probability $\mathbb{P}$ in the filtration $\mathbb{F}$, if the non decreasing process $\sqrt{\transp{H}(\centerdot[X,\transp X]){H}}$ is $(\mathbb{P},\mathbb{F})$ locally integrable. For such an integrable process ${H}$, the stochastic integral $\transp{H}\centerdot X$ is well-defined and the bracket process of $\transp{H}\centerdot X$ can be computed using \cite[Remarque(4.36) and Proposition(4.68)]{Jacodlivre}. Note that two different predictable processes may produce the same stochastic integral with respect to $X$. In this case, we say that they are in the same equivalent class.

Again another notion of stochastic integral is needed, i.e., the stochastic integral with respect to a compensated integer valued random measure $\mu-\nu$. We refer to \cite{HWY, Jacodlivre,JacShi} for its definition and the fundamental properties. In particular, we denote by $\mathscr{G}(\mathbb{F},\mu)$ the space of $(\mu-\nu)$ $_*$-integrable $\mathbb{F}$ predictable functions. To distinguish the different type of stochastic integrals, we mention the stochastic integral with respect to a compensated integer valued random measure $\mu-\nu$ as stochastic $_*$-integral, whilst the stochastic integral with respect to a semimartingale will be mentioned as stochastic $\centerdot$-integral.

\textbf{Caution.}
Note that some same notations are used in different parts of the paper for different meaning, especially the notations $X,Y,H,G,S,T$, $\mu$ or $\nu$.

\

\section{Finite predictable constraint}\label{sectionFPCC}

This section is devoted to the condition of finite predictable constraint for integer valued random measures. Recall that $\mathbb{F}$ is a filtration on the probability space $(\Omega,\mathcal{A},\mathbb{P})$, satisfying the usual conditions.

\subsection{Martingale projection property}\label{pjp}

When we compute the predictable bracket $[Y,M]^{\mathbb{F}\cdot p}$ for two local martingales $M,Y$ ($M$ being, say, locally bounded), we may need to substitute $Y$ by its \texttt{"}orthogonal\texttt{"} projection onto the stable space generated by $M$ : $$
[Y,M]^{\mathbb{F}\cdot p}=[H\centerdot M,M]^{\mathbb{F}\cdot p},
$$
with a stochastic integral $H\centerdot M$ with respect to $M$. It is however not always possible, as explained in \cite{anselstricker}. On the other hand, as a consequence of \cite[Theorem (3.75)]{Jacodlivre}, we have a general projection formula for stochastic $_*$-integrals.

\bl\label{projectionlemma}
Let $M$ be a multiple dimensional purely discontinuous $(\mathbb{P},\mathbb{F})$ local martingale. Let $\mu$ be its jump measure with $(\mathbb{P},\mathbb{F})$ compensator $\nu$. For any real $(\mathbb{P},\mathbb{F})$ local martingale $Y$ such that $[Y,M]$ is $(\mathbb{P},\mathbb{F})$ locally integrable. There exists a $g\in\mathscr{G}(\mathbb{F},\mu)$ such that $[g{_*}(\mu-\nu), M]$ is locally integrable and$$
[Y,M]^{\mathbb{F}\cdot p}=[g{_*}(\mu-\nu), M]^{\mathbb{F}\cdot p}.
$$  
\el

\textbf{Proof.}
Denote by $\mathsf{M}$ the Dolean-Dade measure associated with $\mu$. Let $(U_n)_{n\in\mathbb{N}}$ be a sequence of $\mathbb{F}$ stopping times, tending to the infinity, such that $\mathbb{E}[\int_0^{U_n}|\mathsf{d}[Y,M_h]|]<\infty$ for every component $M_h$ of $M$ and for every $n\in\mathbb{N}$. This implies$$
\mathsf{M}[|\Delta Y| |x_h|\ind_{[0,U_n]}]
=
\mathbb{E}[\sum_{0<s\leq U_n}|\Delta_s Y\ \Delta_s M_h|\ind_{\{\Delta_sM\neq 0\}}]
=
\mathbb{E}[\int_0^{U_n}|\mathsf{d}[Y,M_h]|]<\infty.
$$
With the notations in \cite[Theorem (3.75)]{Jacodlivre} let $$
U = \mathsf{M}[\Delta Y |\widetilde{\mathcal{P}}],\
g  = U + \frac{\widehat{U}}{1-a},\
V = \Delta Y - U.
$$
Then, $g\in\mathscr{G}(\mathbb{F},\mu)$ and$$
Y = g{_*}(\mu-\nu) + V{_*}\mu + Y',
$$
where $Y'$ is a local martingale pathwisely orthogonal to $M$, i.e. $[Y',M]\equiv 0$. Consider $[V{_*}\mu,M_h]$. We have$$
\dcb
&&\mathbb{E}[\int_0^{U_n}|\mathsf{d}[V{_*}\mu,M_h]|]
=
\mathbb{E}[\sum_{0<s\leq U_n}|(\Delta_s Y-U(s,\Delta_s M)) \Delta_s M_h|\ind_{\{\Delta_sM\neq 0\}}]\\

&\leq&
\mathbb{E}[\sum_{0<s\leq U_n}|\Delta_s Y\Delta_s M_h|\ind_{\{\Delta_sM\neq 0\}}]
+
\mathbb{E}[\sum_{0<s\leq U_n}|U(s,\Delta_s M) \Delta_s M_h|\ind_{\{\Delta_sM\neq 0\}}]\\

&=&
\mathsf{M}[|\Delta Y| |x_h|\ind_{[0,U_n]}]+\mathsf{M}[|U| |x_h|\ind_{[0,U_n]}]\\

&\leq&
2\mathsf{M}[|\Delta Y| |x_h|\ind_{[0,U_n]}]
<\infty.
\dce
$$
This means that $[V{_*}\mu,M_h]^{\mathbb{F}\cdot p}$ is defined. But for any $\mathbb{F}$ stopping time $S$, $$
\mathbb{E}[[V{_*}\mu,M_h]_{S\wedge U_n}]
=
\mathsf{M}[(\Delta Y-U) x_h \ind_{[0,S\wedge U_n]}]=0,
$$
i.e. $[V{_*}\mu,M_h]^{\mathbb{F}\cdot p}=0$. As $[Y,M_h]^{\mathbb{F}\cdot p}$ exists, necessarily $[g{_*}(\mu-\nu),M_h]^{\mathbb{F}\cdot p}$ exists and 
$$
[Y,M]^{\mathbb{F}\cdot p}=[g{_*}(\mu-\nu), M]^{\mathbb{F}\cdot p}.\ \ok
$$

In the light of Lemma \ref{projectionlemma}, we understand that, to have the martingale projection property for stochastic ${\centerdot}$-integrals, it is enough to find conditions which make the stochastic $_*$-integrals $g{_*}(\mu-\nu)$ to become stochastic ${\centerdot}$-integrals $H\centerdot M$.

\

\subsection{The definition}\label{irm}

We now study the general problem of the transformation from stochastic $_{*}$-integrals into stochastic ${\centerdot}$-integrals. We recall the basic vocabulary about integer valued random measures. Let $\mathtt{E}$ be an Euclidean space. An $\mathbb{F}$ optional random measure $\mu$ on $\mathbb{R}_+\times\mathtt{E}$ is said to be integer valued (cf. \cite{HWY,Jacodlivre}), if there exists an $\mathbb{F}$ optional thin set $\mathtt{D}$ (the time support set) and an $\mathtt{E}$-valued $\mathbb{F}$ optional process $\beta$ (the space location process) such that$$
\mu[\mathtt{A}] = \sum_{s>0}\ind_{\{(s,\beta_s)\in\mathtt{A}\}}\ind_{\{s\in\mathtt{D}\}},\ \forall \mathtt{A}\in\mathcal{B}(\mathbb{R}_+\times\mathtt{E}).
$$ 
We make use of the results in \cite[Chapiter XI section 1]{HWY}, also in \cite[Chapiter II section 1]{JacShi}. In this paper, the integer valued random measures $\mu$ are always supposed to be $\sigma$-finite on the predictable $\sigma$-algebra and to have an $\mathbb{F}$ compensator $\nu$ satisfying
\begin{equation}\label{rmcond}
\nu[\{0\}\times\mathtt{E}]
=
\nu[\mathbb{R}_+\times\{\boldsymbol{0}\}]=0,\
(|x|^2\wedge 1){_*}\nu_t <\infty, \ t\in\mathbb{R}_+.
\end{equation} 
Note that these conditions are satisfied by the jump measure of any semimartingale (cf. \cite[Chapter II, Proposition 2.9]{JacShi}). Recall that $\mathscr{G}(\mathbb{F},\mu)$ denotes the $(\mu-\nu)$ ${_*}$-integrable predictable functions.

Here is the notion which makes the transformation from stochastic $_{*}$-integrals into stochastic ${\centerdot}$-integrals possible.

\bd\label{FPCC}
We say that an integer valued $\mathbb{F}$ optional random measure $\mu$ satisfies the finite predictable constraint condition, if the space location process $\beta$ is confined in a finite $\mathbb{F}$ predictable constraint, i.e., if there exist a finite number (say $\mathsf{n}$) of $\mathtt{E}$-valued $\mathbb{F}$ predictable processes $\alpha_k, 1\leq k\leq \mathsf{n}$, such that, at any time, the value of $\beta$ coincides with one of the values $\alpha_{k}$ or $\boldsymbol{0}$:$$
\beta\in\{\boldsymbol{0},\alpha_1,\ldots,\alpha_\mathsf{n}\}.
$$
\ed

\brem
Note that, in the case of a finite predictable constraint for $\beta$, we can modify the $\mathtt{E}$-valued $\mathbb{F}$ predictable constraint processes $\alpha_k, 1\leq k\leq \mathsf{n}$, to write 
\begin{equation}\label{beta=alpha}
\beta
= 
\sum_{k=1}^{\mathsf{n}}\alpha_k\ind_{\{\beta=\alpha_k\}} + \boldsymbol{0}\ind_{\{\beta\neq \alpha_{k}, \forall k\}}.
\end{equation}
We accept some set $\{\beta=\alpha_k\}$ empty.
\erem

\

\subsection{The main result}

We consider an integer valued random measure $\mu$ with its time support $\mathtt{D}$ and its space location process $\beta$ and its $\mathbb{F}$ compensator $\nu$. Suppose the finite predictable constraint condition with the constraint processes $\alpha_k, 1\leq k\leq \mathsf{n}.$ Let $e_k$ ($1\leq k\leq \mathsf{n}$) be a bounded continuous real function such that $e_k(\alpha_k)\neq 0$ on the time support set $\mathtt{D}$ of $\mu$, and $|e_k(x)|\leq c (|x|\wedge 1), x\in\mathtt{E},$ for some constant $c$.

\bethe\label{specif}
Suppose that $\mu$ satisfies the finite $\mathbb{F}$ predictable constraint condition with constraint processes $\alpha_k, 1\leq k\leq \mathsf{n},$ satisfying identity (\ref{beta=alpha}). Suppose $\mathtt{D}\subset\{\beta\neq \boldsymbol{0}\}$.  For $1\leq k\leq \mathsf{n}$, let $e_k$ be the above defined functions. 
\ebe
\item
Let $u_k(s,x)=e_k(x)\ind_{\{x=\alpha_{k,s}\}}, s\geq 0, x\in\mathtt{E},$ and $X_k=u_k{_*}(\mu-\nu), 1\leq k\leq \mathsf{n}$. Then, $X_k$ are well-defined locally bounded $\mathbb{F}$ local martingales.
\item
For an element $g\in \mathscr{G}(\mathbb{F},\mu)$, let $\frac{g(\cdot,\alpha)}{\mathsf{e}(\alpha)}\ind_{\{e(\alpha)\neq 0\}}$ denote the vector valued process composed of $g(\cdot,\alpha_{k})\frac{1}{e_k(\alpha_k)}\ind_{\{e_k(\alpha_k)\neq 0\}}, 1\leq k\leq \mathsf{n}$. Then, $\frac{g(\cdot,\alpha)}{\mathsf{e}(\alpha)}\ind_{\{e(\alpha)\neq 0\}}$ is integrable with respect to the vector valued process $X=(X_k)_{1\leq k\leq \mathsf{n}}$.
\item
For any element $g\in \mathscr{G}(\mathbb{F},\mu)$,
$$
g{_*}(\mu-\nu)
=
\transp (\frac{g(\cdot,\alpha)}{\mathsf{e}(\alpha)}\ind_{\{e(\alpha)\neq 0\}})\centerdot X.
$$
\dbe 
We have the identity $$
\{g{_*}(\mu-\nu): g\in \mathscr{G}(\mathbb{F},\mu)\}
=
\{\transp H\centerdot X: \mbox{ $H$ is $X$-integrable}\}.
$$
\ethe

\textbf{Proof.} Using the notations in \cite[Definition 11.16]{HWY}, we must prove that the process $(\sqrt{\sum_{s\leq t}\widetilde{u}_{k,s}^2}, t\geq 0)$ is locally integrable, where $$
\widetilde{u}_{k} = u_{k}(\cdot, \beta)\ind_{\mathtt{D}} - \widehat{u}_{k}.
$$
We consider separately $\sum_{s\leq t, s\in\mathtt{D}}{u}_{k}(s,\beta_{s})^2$ and $\sum_{s\leq t}\widehat{u}_{k,s}^2$. For any stopping time $T$ such that $
\mathbb{E}[(|x|^2\wedge 1){_{*}}\nu_{T}]<\infty,
$
we have$$
\mathbb{E}[\sum_{s\leq T, s\in\mathtt{D}}{u}_{k}(s,\beta_{s})^2]
=
\mathbb{E}[{u}_{k}^2{_{*}}\mu_{T}]
\leq
c^2\mathbb{E}[(|x|^2\wedge 1){_{*}}\mu_{T}]
=
c^2\mathbb{E}[(x^2\wedge 1){_{*}}\nu_{T}]<\infty.
$$
On the other hand, we know that $\{\widehat{u}_{k}\neq 0\}$ is a predictable thin set $\mathtt{J}$ (cf. \cite[Theorem 11.14]{HWY}). Hence,
$$
\mathbb{E}[\sum_{s\leq T}\widehat{u}_{k,s}^2]
=
\mathbb{E}[\sum_{s\leq T, s\in \mathtt{J}}\widehat{u}_{k,s}^2]
\leq 
\mathbb{E}[\sum_{s\leq T, s\in \mathtt{J}}\int_{[s]\times \mathtt{E}} {u}_{k}^2 \mathsf{d}\nu]
=
\mathbb{E}[\int_{[0,T]\times \mathtt{E}} \ind_{J}{u}_{k}^2 \mathsf{d}\nu]
\leq
c^2\mathbb{E}[\int_{[0,T]\times \mathtt{E}} (|x|^2\wedge 1) \mathsf{d}\nu]<\infty.
$$
Because of the conditions in (\ref{rmcond}), we conclude that the process $(\sqrt{\sum_{s\leq t}\widetilde{u}_{k,s}^2}, t\geq 0)$ is locally integrable and the local martingale $X_k$ is well defined. 

With identity (\ref{beta=alpha}), necessarily $\alpha_{k'}\neq \alpha_{k}$ on $\{\beta\neq \boldsymbol{0}, \beta=\alpha_{k}\}$, for all $k'\neq k$. As $\mathtt{D}\subset\{\beta\neq \boldsymbol{0}\}$, the sets $\{s\in\mathtt{D},\beta_s=\alpha_k\}, 1\leq k\leq \mathsf{n}$, are mutually disjoint, and $\{s\in\mathtt{D},\forall k, \beta_s\neq\alpha_k\}=\emptyset$. So, for a $g\in \mathscr{G}(\mathbb{F},\mu)$, for any $\mathbb{F}$ stopping time $T$, $$
g(T,\beta_T)\ind_{\{T\in\mathtt{D}\}}
=
\sum_{k=1}^{\mathsf{n}}g(T,\alpha_{k,T})\ind_{\{\beta_T=\alpha_{k,T}\}}\ind_{\{T\in\mathtt{D}\}}.
$$
We compute the jump at an $\mathbb{F}$ totally inaccessible stopping time $T$ on $\{T<\infty\}$.$$
\dcb
&&\Delta_T (g{_*}(\mu-\nu))
=
g(T,\beta_T)\ind_{\{T\in\mathtt{D}\}}
=
\sum_{k=1}^{\mathsf{n}}g(T,\alpha_{k,T})\ind_{\{\beta_T=\alpha_{k,T}\}}\ind_{\{T\in\mathtt{D}\}}\\
&=&
\sum_{k=1}^{\mathsf{n}}g(T,\alpha_{k,T})\frac{1}{e_k(\alpha_{k,T})}\ind_{\{e_k(\alpha_{k,T})\neq 0\}} \Delta_T X_k.
\dce
$$
We compute next the jump at an $\mathbb{F}$ predictable stopping time $T$ on $\{T<\infty\}$.$$
\dcb
&&\Delta_T (g{_*}(\mu-\nu))\\
&=&
g(T,\beta_T)\ind_{\{T\in\mathtt{D}\}} - \mathbb{E}[g(T,\beta_T)\ind_{\{T\in\mathtt{D}\}}|\mathcal{F}_{T-}]\\
&=&
\sum_{k=1}^{\mathsf{n}}g(T,\alpha_{k,T})\ind_{\{\beta_T=\alpha_{k,T}\}}\ind_{\{T\in\mathtt{D}\}}
- \mathbb{E}[\sum_{k=1}^{\mathsf{n}}g(T,\alpha_{k,T})\ind_{\{\beta_T=\alpha_{k,T}\}}\ind_{\{T\in\mathtt{D}\}}|\mathcal{F}_{T-}]\\

&=&
\sum_{k=1}^{\mathsf{n}}g(T,\alpha_{k,T})(\ind_{\{\beta_T=\alpha_{k,T}\}}\ind_{\{T\in\mathtt{D}\}}
-
\mathbb{E}[\ind_{\{\beta_T=\alpha_{k,T}\}}\ind_{\{T\in\mathtt{D}\}}|\mathcal{F}_{T-}])\\

&=&
\sum_{k=1}^{\mathsf{n}}g(T,\alpha_{k,T})\frac{1}{e_k(\alpha_{k,T})}\ (e_k(\alpha_{k,T})\ind_{\{\beta_T=\alpha_{k,T}\}}\ind_{\{T\in\mathtt{D}\}}
-
\mathbb{E}[e_k(\alpha_{k,T})\ind_{\{\beta_T=\alpha_{k,T}\}}\ind_{\{T\in\mathtt{D}\}}|\mathcal{F}_{T-}])\\

&=&
\sum_{k=1}^{\mathsf{n}}g(T,\alpha_{k,T})\frac{1}{e_k(\alpha_{k,T})}\ (e_k(\beta_{T})\ind_{\{\beta_T=\alpha_{k,T}\}}\ind_{\{T\in\mathtt{D}\}}
-
\mathbb{E}[e_k(\beta_{T})\ind_{\{\beta_T=\alpha_{k,T}\}}\ind_{\{T\in\mathtt{D}\}}|\mathcal{F}_{T-}])\\

&=&
\sum_{k=1}^{\mathsf{n}}g(T,\alpha_{k,T})\frac{1}{e_k(\alpha_{k,T})}\ind_{\{e_k(\alpha_{k,T})\neq 0\}} \Delta_T X_k.
\dce
$$
We obtain the identity$$
\Delta (g{_*}(\mu-\nu))
=
\sum_{k=1}^{\mathsf{n}}g(\cdot,\alpha_{k})\frac{1}{e_k(\alpha_{k})}\ind_{\{e_k(\alpha_{k})\neq 0\}} \Delta X_k.
$$
This identity shows firstly that the process $\frac{g(\cdot,\alpha)}{\mathsf{e}(\alpha_{})}\ind_{\{e(\alpha)\neq 0\}}$ is $X$-integrable. Secondly, by \cite[Theorem 7.23]{HWY}, we have the equality$$
g{_*}(\mu-\nu)
=
\transp \left(\frac{g(\cdot,\alpha)}{\mathsf{e}(\alpha_{})}\ind_{\{e(\alpha)\neq 0\}}\right)\centerdot X.\ \ok
$$

\

\subsection{Case of random measure with accessible time support}\label{pjpacc}

To represent the space of stochastic $_*$-integrals $\{g{_*}(\mu-\nu): g\in\mathscr{G}(\mu)\}$ with stochastic $\centerdot$-integrals under the finite predictable constraint condition, we may employ different \texttt{"}representation\texttt{"} local martingales than that defined in Theorem \ref{specif}, especially when $\mu$ is the jump measure of a local martingale.

Consider an integer valued random measure $\mu$ with compensator $\nu$. We now cut the process $\beta$ into pieces in time and in space. Suppose that there exist a sequence of mutually avoiding $\mathbb{F}$ predictable stopping times $(T_n)_{1\leq n<\mathsf{N}}$ ($\mathsf{N}\leq \infty$) such that the time support set is given by $\mathtt{D}=( \cup_{1\leq n<\mathsf{N}}[T_n])\cap \{\beta\neq \boldsymbol{0}\}$. We consider the family of random variables $(\beta_{T_n}, 1\leq n<\mathsf{N})$. Suppose in addition that there exist a positive integer $\mathsf{n}$ such that, for every $1\leq n<\mathsf{N}$, there exist $\mathtt{E}$-valued $\mathcal{F}_{T_n-}$ measurable $\alpha_{n,k}, 1\leq k\leq \mathsf{n}$, and an $\mathcal{F}_{T_n}$ measurable partition $(A_{n,1},\ldots,A_{n,\mathsf{n}})$ (possibly some empty sets) such that $\beta_{T_n}$ is cut into
\begin{equation}\label{betapartition}
\beta_{T_n}
= 
\sum_{k=1}^{\mathsf{n}}\alpha_{n,k}\ind_{A_{n,k}}.
\end{equation}
Let $(a_n)_{1\leq n<\mathsf{N}}$ be any series of non vanishing random variables such that $a_n\in\mathcal{F}_{T_n-}$ and the following expression 
$$
Y_k=\sum_{n=1}^{\mathsf{N}-}a_n(\ind_{A_{n,k}}\ind_{[T_n,\infty)}-(\ind_{A_{n,k}}\ind_{[T_n,\infty)})^{\mathbb{F}\cdot p}),\ 1\leq k\leq \mathsf{n},
$$
defines an $\mathsf{n}$-dimensional $\mathbb{F}$ local martingale $Y$. Set $\alpha_k=\sum_{n=1}^{\mathsf{N}-}\alpha_{n,k}\ind_{[T_n]}, 1\leq k\leq \mathsf{n},$ and $
G=\sum_{n=1}^{\mathsf{N}-}\frac{1}{a_n}\ind_{[T_n]}.
$
For any $g\in\mathscr{G}(\mathbb{F},\mu)$, denote by $g(\cdot,\alpha)\ind_{\{\alpha\neq 0\}}$ the vector valued process $(g(\cdot,\alpha_{k})\ind_{\{\alpha_{k}\neq \boldsymbol{0}\}})_{1\leq k\leq \mathsf{n}}$. 

\bethe\label{specif2}
Under the above conditions, for any $g\in\mathscr{G}(\mathbb{F},\mu)$, $Gg(\cdot,\alpha)\ind_{\{\alpha\neq 0\}}$ is $Y$-integrable and $$
g{_*}(\mu-\nu)
=
G\ \transp g(\cdot,\alpha)\ind_{\{\alpha\neq 0\}}\centerdot Y.
$$
\ethe

\brem
If $\beta$ satisfies condition (\ref{betapartition}), (a version of) $\beta$ satisfied the finite predictable constraint. \ok
\erem

\textbf{Proof.} 
For $g\in \mathscr{G}(\mathbb{F},\mu)$, as in the proof of Theorem \ref{specif}, we compute the jumps at one of the $\mathbb{F}$ predictable stopping times $T=T_n<\infty$ ($1\leq n<\mathsf{N}$).$$
\dcb
&&\Delta_T (g{_*}(\mu-\nu))\\
&=&
g(T,\beta_T)\ind_{\{T\in\mathtt{D}\}} - \mathbb{E}[g(T,\beta_T)\ind_{\{T\in\mathtt{D}\}}|\mathcal{F}_{T-}]\\

&=&
g(T,\beta_T)\ind_{\{\beta_T\neq \boldsymbol{0}\}} - \mathbb{E}[g(T,\beta_T)\ind_{\{\beta_T\neq \boldsymbol{0}\}}|\mathcal{F}_{T-}]\\

&=&
\sum_{k=1}^{\mathsf{n}}g(T,\alpha_{k,T})\ind_{A_{n,k}}\ind_{\{\alpha_{k,T}\neq \boldsymbol{0}\}}
- \mathbb{E}[\sum_{k=1}^{\mathsf{n}}g(T,\alpha_{k,T})\ind_{A_{n,k}}\ind_{\{\alpha_{k,T}\neq \boldsymbol{0}\}}|\mathcal{F}_{T-}]\\

&=&
\sum_{k=1}^{\mathsf{n}}g(T,\alpha_{k,T})\ind_{\{\alpha_{k,T}\neq \boldsymbol{0}\}}(\ind_{A_{n,k}}
-
\mathbb{E}[\ind_{A_{n,k}}|\mathcal{F}_{T-}])\\

&=&
\sum_{k=1}^{\mathsf{n}}g(T,\alpha_{k,T})\ind_{\{\alpha_{k,T}\neq \boldsymbol{0}\}}\frac{1}{a_n}\ \Delta_T Y_k\\

&=&
\sum_{k=1}^{\mathsf{n}}g(T,\alpha_{k,T})\ind_{\{\alpha_{k,T}\neq \boldsymbol{0}\}}G_{T}\ \Delta_T Y_k.

\dce
$$
From this jump identity, we conclude that $Gg(\cdot,\alpha)\ind_{\{\alpha\neq 0\}}$ is $Y$-integrable, and, by \cite[Theorem 7.23]{HWY},$$
g{_*}(\mu-\nu)
=
G\ \transp (g(\cdot,\alpha)\ind_{\{\alpha\neq 0\}})\centerdot Y.\ \ok
$$ 

\

Suppose, on top of the conditions in Theorem \ref{specif2}, that, for a $d$-dimensional $\mathbb{F}$ purely discontinuous local martingale $M$, $\mathtt{D}=\{\Delta M\neq \boldsymbol{0}\}$ and $\beta=\Delta M$. For a fixed $1\leq n<\mathsf{N}$, let us view the $\mathcal{F}_{T_n-}$ measurable random variables as constants and consider the space $\mathfrak{L}_{n}$ of real functions generated by the indicators $\ind_{A_{n,k}}, 1\leq k\leq \mathsf{n}$. The space $\mathfrak{L}_{n}$ is of finite dimension and there exists a natural linear map from $\mathbb{R}^\mathsf{n}$ ($\mathcal{F}_{T_n-}$ measurable random variables being constant) onto the space $\mathfrak{L}_{n}$. Condition (\ref{betapartition}) says that the components $\Delta_{T_n} M_{i}, 1\leq i\leq d$, of $\Delta_{T_n} M$ are elements in the space $\mathfrak{L}_{n}$:$$
\Delta_{T_n}M_{i}
= 
\sum_{k=1}^{\mathsf{n}}\alpha_{n,i,k}\ind_{A_{n,k}},
$$
where $\alpha_{n,i,k}$ denotes the $i$th component of $\alpha_{n,k}\in\mathbb{R}^d$. Introduce the (vertical) vectors $\gamma_{n,i}$ of the components $(\alpha_{n,i,k})_{1\leq k\leq \mathsf{n}}$ ($\gamma_{n,i}$ being a representative of $\Delta_{T_n} M_{i}$ in $\mathbb{R}^\mathsf{n}$). Denote $p_{n,k}=\mathbb{E}[\ind_{A_{n,k}}|\mathcal{F}_{T_n-}]$ and $p_n$ the (vertical) vectors of the components $(p_{n,k})_{1\leq k\leq \mathsf{n}}$. We have for $1\leq i\leq d$
$$
0=\mathbb{E}[\Delta_{T_n}M_i|\mathcal{F}_{T_n-}]
=
\sum_{k=1}^\mathsf{n}\alpha_{n,i,k}\mathbb{E}[\ind_{A_{n,k}}|\mathcal{F}_{T_n-}]
=
\sum_{k=1}^\mathsf{n}\alpha_{n,i,k}p_{n,k}
=
\transp \gamma_{n,i} p_{n},
$$
i.e., $\gamma_{n,i}$ is orthogonal to $p_n$. Let $P_{n}=\sum_{k=1}^{\mathsf{n}} p_{n,k}\ind_{A_{n,k}}$.

\bethe
Under the above condition, suppose that the vectors $\gamma_{n,i}, 1\leq i \leq d,$ together with $p_n$ span the whole space $\mathbb{R}^\mathsf{n}$ (so that $\Delta_{T_{n}}M_{i}, 1\leq i\leq d,$ and $P_{n}$ generates $\mathfrak{L}_{n}$ and therefore $\mathsf{n}\leq d+1$). Then, for any $g\in\mathscr{G}(\mathbb{F},\mu)$, there exists a matrix valued $\mathbb{F}$ predictable process $K$ such that $Kg(\cdot,\alpha)\ind_{\{\alpha\neq 0\}}$ is $M$-integrable and
$$
g{_*}(\mu-\nu)
=
\transp (Kg(\cdot,\alpha)\ind_{\{\alpha\neq 0\}})\centerdot M.
$$ 
We have the identity$$
\{g{_*}(\mu-\nu): g\in\mathscr{G}(\mathbb{F},\mu)\}
=
\{\transp H\centerdot M: \mbox{ $H$ is $M$-integrable}\}.
$$
\ethe

\textbf{Proof.}
For $g\in \mathscr{G}(\mathbb{F},\mu)$, as in the preceding proof, we compute the jumps at one of the $\mathbb{F}$ predictable stopping times $T_n<\infty$ ($1\leq n<\mathsf{N}$).$$
\dcb
&&\Delta_{T_n} (g{_*}(\mu-\nu))
=
\sum_{h=1}^{\mathsf{n}}g({T_n},\alpha_{h,{T_n}})\ind_{\{\alpha_{h,{T_n}}\neq 0\}}(\ind_{A_{n,h}}
-
\mathbb{E}[\ind_{A_{n,h}}|\mathcal{F}_{{T_n}-}]).

\dce
$$
Note that$$
(\ind_{A_{n,h}}
-
\mathbb{E}[\ind_{A_{n,h}}|\mathcal{F}_{{T_n}-}])
=
(\ind_{A_{n,h}}
-
p_{n,h})
=
\sum_{k=1}^\mathsf{n}(\delta_{h,k}-p_{n,h})\ind_{A_{n,k}}.
$$
Taking the conditioning with respect to $\mathcal{F}_{{T_n}-}$, we see $$
\sum_{k=1}^\mathsf{n}(\delta_{h,k}-p_{n,h})p_{{n,k}}=0,
$$
i.e., the vector of components $(\delta_{h,k}-p_{n,h})_{1\leq k\leq \mathsf{n}}$ is orthogonal to $p_n$ so that there exists a $\mathcal{F}_{{T_n}-}$-measurable vector $K_{n,h}=(K_{n,i,h})_{1\leq i\leq d}$ such that $$
(\delta_{h,k}-p_{n,h})_{1\leq k\leq \mathsf{n}}
=
\sum_{i=1}^d \gamma_{n,i}K_{n,i,h},
$$
or in other words, the image in $\mathfrak{L}_{n}$ of the vector $(\delta_{h,k}-p_{n,h})_{1\leq k\leq \mathsf{n}}$ is a combination of the $\Delta_{U_{n}}M_{i}$:  
$$
\dcb
&&(\ind_{A_{n,h}}
-
\mathbb{E}[\ind_{A_{n,h}}|\mathcal{F}_{{T_n}-}])
=
\sum_{k=1}^\mathsf{n}(\delta_{h,k}-p_{n,h})\ind_{A_{n,k}}\\

&=&
\sum_{k=1}^\mathsf{n}(\sum_{i=1}^d \gamma_{n,i}K_{n,i,h})_k\ind_{A_{n,k}}
=
\sum_{k=1}^\mathsf{n}\sum_{i=1}^d \alpha_{n,i,k}K_{n,i,h}\ind_{A_{n,k}}\\

&=&
\sum_{i=1}^d K_{n,i,h}\sum_{k=1}^\mathsf{n}\alpha_{n,i,k}\ind_{A_{n,k}}
=
\sum_{i=1}^d K_{n,i,h}\Delta_{T_n}M_i
=
\transp K_{n,h} \Delta_{T_n}M.
\dce
$$
Hence,
$$
\dcb
\Delta_{T_n} (g{_*}(\mu-\nu))
&=&
\sum_{h=1}^{\mathsf{n}}g({T_n},\alpha_{h,{T_n}})\ind_{\{\alpha_{h,{T_n}}\neq 0\}}(\ind_{A_{n,h}}
-
\mathbb{E}[\ind_{A_{n,h}}|\mathcal{F}_{{T_n}-}])\\

&=&
\sum_{h=1}^{\mathsf{n}}g({T_n},\alpha_{h,{T_n}})\ind_{\{\alpha_{h,{T_n}}\neq 0\}}\transp K_{n,h} \Delta_{T_n}M.
\dce
$$
Set $K_n$ to be the matrix $(K_{n,i,h})_{1\leq i\leq d,1\leq h\leq\mathsf{n}}$ and $K=\sum_{n=1}^{\mathsf{N}-}K_{n}\ind_{[T_n]}$. The above jump identity implies that $K g(\cdot,\alpha)\ind_{\{\alpha\neq 0\}}$ is $M$-integrable, and, by \cite[Theorem 7.23]{HWY},$$
g{_*}(\mu-\nu)
=
\transp (Kg(\cdot,\alpha)\ind_{\{\alpha\neq 0\}})\centerdot M.
$$ 
To finish the proof, we recall that $M_{h}= x_{h}{_*}(\mu-\nu)$. \ \ok

\

\subsection{Case of totally inaccessible support}\label{pjpinacc}

Consider always an integer valued random measure $\mu$ with its compensator $\nu$, satisfying the finite predictable constraint condition with constraint processes $\alpha_k, 1\leq k\leq \mathsf{n},$ satisfying identity (\ref{beta=alpha}). Suppose that $\mathtt{D}=\{\beta\neq \boldsymbol{0}\}=\cup_{1\leq n<\mathsf{N}}[S_n]$, where $S_n$ are mutually avoiding $\mathbb{F}$ totally inaccessible stopping times and $\mathsf{N}$ an finite or infinite integer.

Suppose in addition that, for a $d$-dimensional $\mathbb{F}$ purely discontinuous local martingale $M$, $\mathtt{D}=\{\Delta M\neq \boldsymbol{0}\}$ and $\beta=\Delta M$. ($M$ is then quasi-left continuous. See \cite[Theorem 4.23]{HWY}). The processes $\alpha_{k}$ are therefore $d$-dimensional vectors. Let $\alpha_{i,k}$ denote the $i$th component of $\alpha_{k}$ for $1\leq i\leq d$. As in the preceding paragraph, we define the vector $\gamma_i, 1\leq i\leq d,$ to be the vector of the components $(\alpha_{i,k})_{1\leq k\leq \mathsf{n}}$. Denote by $\gamma$ the matrix of columns $\gamma_i$'s.

\bethe\label{specif3}
Under the above conditions, suppose that, for any $S_n$, the vectors $\gamma_{i,S_n}, 1\leq i\leq d$, span the whole space $\mathbb{R}^\mathsf{n}$ (so that $\mathsf{n}\leq d$). Then, for any $g\in\mathscr{G}(\mathbb{F},\mu)$, there exists a matrix valued $\mathbb{F}$ predictable process $K$ such that $Kg(\cdot,\alpha)\ind_{\{\alpha\neq 0\}}$ is $M$-integrable and
$$
g{_*}(\mu-\nu)
=
\transp (Kg(\cdot,\alpha)\ind_{\{\alpha\neq 0\}})\centerdot M.\ \ok
$$ 
We have the identity$$
\{g{_*}(\mu-\nu): g\in\mathscr{G}(\mathbb{F},\mu)\}
=
\{\transp H\centerdot M: \mbox{ $H$ is $M$-integrable}\}.
$$ 
\ethe

\textbf{Proof.}
Consider the canonical basis $(\epsilon_1,\ldots,\epsilon_\mathsf{n})$ in $\mathbb{R}^\mathsf{n}$. Note that the set $$
\mathtt{A}=\{\mbox{$\gamma_{i}, 1\leq i\leq d$, span the whole space $\mathbb{R}^\mathsf{n}$}\}
$$
is $\mathbb{F}$ predictable. There exists an $\mathbb{F}$ predictable $d\times\mathsf{n}$-matrix valued process $K$ such that $$
(\epsilon_1,\ldots,\epsilon_\mathsf{n})
=
\gamma K\ \mbox{ on $\mathtt{A}$}.
$$
For any element $g\in \mathscr{G}(\mathbb{F},\mu)$,
we compute the jump at a stopping time $S_n<\infty$.$$
\dcb
&&\Delta_{S_n} (g{_*}(\mu-\nu))
=
g({S_n},\beta_{S_n})
=
\sum_{h=1}^{\mathsf{n}}g({S_n},\alpha_{h,{S_n}})\ind_{\{\alpha_{h,{S_n}}\neq \boldsymbol{0}\}}\ind_{\{\beta_{S_n}=\alpha_{h,{S_n}}\}}.
\dce
$$
Note (as in the preceding proof) that$$
\dcb
&&\ind_{\{\beta_{S_n}=\alpha_{h,{S_n}}\}}
=
\sum_{k=1}^\mathsf{n}\delta_{h,k}\ind_{\{\beta_{S_n}=\alpha_{k,{S_n}}\}}
=
\sum_{k=1}^\mathsf{n}\sum_{i=1}^d \gamma_{i,k,S_n}K_{i,h,S_n}\ind_{\{\beta_{S_n}=\alpha_{k,{S_n}}\}}\\

&=&
\sum_{k=1}^\mathsf{n}\sum_{i=1}^d \alpha_{i,k,S_n}K_{i,h,S_n}\ind_{\{\beta_{S_n}=\alpha_{k,{S_n}}\}}

=
\sum_{i=1}^dK_{i,h,S_n}\sum_{k=1}^\mathsf{n} \alpha_{i,k,S_n}\ind_{\{\beta_{S_n}=\alpha_{k,{S_n}}\}}\\

&=&
\sum_{i=1}^dK_{i,h,S_n}\Delta_{S_n}M_i.

\dce
$$
We conclude
$$
\dcb
&&\Delta_{S_n} (g{_*}(\mu-\nu))
=
\sum_{h=1}^{\mathsf{n}}g({S_n},\alpha_{h,{S_n}})\ind_{\{\alpha_{h,{S_n}}\neq \boldsymbol{0}\}}\sum_{i=1}^dK_{i,h,S_n}\Delta_{S_n}M_i\\
&=&
\transp (g(\cdot,\alpha)\ind_{\{\alpha\neq 0\}})_{S_n}\transp K_{S_n} \Delta_{S_n}M.
\dce
$$
This jump identity implies that $K g(\cdot,\alpha)\ind_{\{\alpha\neq 0\}}$ is $M$-integrable, and, by \cite[Theorem 7.23]{HWY},$$
g{_*}(\mu-\nu)
=
\transp (Kg(\cdot,\alpha)\ind_{\{\alpha\neq 0\}})\centerdot M.\ \ok
$$

\

\section{Martingale representation property}\label{mrp}

Fix a stochastic basis $(\Omega,\mathbb{F},\mathbb{P})$.  

\bd\label{assump-mrt}
We say that the martingale representation property holds in the filtration $\mathbb{F}$ (under the probability $\mathbb{P}$) with respect to a $d$-dimensional representation process $W$, if $W$ is a $(\mathbb{P},\mathbb{F})$ local martingale, and all $(\mathbb{P},\mathbb{F})$ local martingale is a stochastic integral with respect to $W$. We also say that $W$ possesses the martingale representation property in $\mathbb{F}$. We say simply that the martingale representation property holds in $\mathbb{F}$, if there exists some $(\mathbb{P},\mathbb{F})$ local martingale which possesses the martingale representation property in $\mathbb{F}$.
\ed

\

\subsection{Conditional multiplicity}\label{multiplicity}

The martingale representation property imposes finite conditional multiplicity of $\mathcal{F}_{R}$ with respect to $\mathcal{F}_{R-}$ (a notion introduced in \cite[section 3]{BEKSY} to quantify the randomness of $\mathcal{F}_{R}$ when $\mathcal{F}_{R-}$ is given).

\bl\label{generating-family}
Suppose the martingale representation property in $\mathbb{F}$ with a $d$-dimensional representation process $W$.
Let $R$ be a $\mathbb{F}$ stopping time. Consider the random variables in $\mathcal{F}_{R-}$ as constants. If $R$ is predictable, the family of random variables $\Delta_R W_h, 1\leq h\leq d$, generates on $\{R<\infty\}$ (modulo $\mathcal{F}_{R-}$) all integrable random variables $\xi$ in $\mathcal{F}_{R}$ whose conditional expectation $\mathbb{E}[\xi|\mathcal{F}_{R-}]=0$. If $R$ is totally inaccessible, the family of $\Delta_R W_h, 1\leq h\leq d$ generates on $\{R<\infty\}$ (modulo $\mathcal{F}_{R-}$) all integrable random variables $\xi$ in $\mathcal{F}_{R}$.
\el

\textbf{Proof.} For any integrable $\xi\in\mathcal{F}_R$, the process $\xi\ind_{[R,\infty)}-(\xi\ind_{[R,\infty)})^{\mathbb{F}\cdot p}$ is a martingale. By martingale representation property, there exists an $\mathbb{F}$-predictable process $H$ such that $\xi\ind_{[T,\infty)}-(\xi\ind_{[T,\infty)})^{\mathbb{F}\cdot p}=\transp H\centerdot W$. Therefore, $$
\xi = \sum_{h=1}^{d}(H_{R})_h\Delta_R W_h+\Delta_R (\xi\ind_{[R,\infty)})^{\mathbb{F}\cdot p}
$$ 
on $\{R<\infty\}$. If $R$ is predictable and $\mathbb{E}[\xi|\mathcal{F}_{R-}]=0$, the process $(\xi\ind_{[R,\infty)})^{\mathbb{F}\cdot p}=0$. If $R$ is totally inaccessible, $\Delta_R (\xi\ind_{[R,\infty)})^{\mathbb{F}\cdot p}=0$. The lemma is proved. \ok

\bl\label{partition}
Suppose the martingale representation property in $\mathbb{F}$ with a $d$-dimensional representation process $W$.
If $R$ is $\mathbb{F}$ predictable, there exists a partition $(A_0,A_1,A_2,\ldots,A_{d})$ (where some $A_i$ may be empty) such that $$
\mathcal{F}_R=\mathcal{F}_{R-}\vee\sigma(A_0,A_1,A_2,\ldots,A_{d}),
$$
i.e. the conditional multiplicity of $\mathcal{F}_R$ with respect to $\mathcal{F}_{R-}$ is equal to or smaller then $d+1$. If $R$ is $(\mathbb{P},\mathbb{F})$ totally inaccessible, there exists a partition $(B_1,B_2,\ldots,B_{d})$ (where some $B_j$ may be empty) such that $$
\mathcal{F}_R=\mathcal{F}_{R-}\vee\sigma(B_1,B_2,\ldots,B_{d}),
$$
i.e. the conditional multiplicity of $\mathcal{F}_R$ with respect to $\mathcal{F}_{R-}$ is equal to or smaller then $d$. 
\el

\textbf{Proof.} Consider the case of a predictable $R$. Because of Lemma \ref{generating-family}, we can apply \cite[Proposition 12]{BEKSY} to have a partition $(A'_0,A'_1,A'_2,\ldots,A'_{d})$ of $\{R<\infty\}$ such that $$
\{R<\infty\}\cap\mathcal{F}_R=\{R<\infty\}\cap(\mathcal{F}_{R-}\vee\sigma(A'_0,A'_1,A'_2,\ldots,A'_{d})).
$$ 
Since $\{R=\infty\}\cap\mathcal{F}_R=\{R=\infty\}\cap\mathcal{F}_{R-}$, the lemma is verified, if we take $A_i=A'_i$ for $0\leq i<d$ and $A_d=A'_d\cup\{R=\infty\}$.

The case of a totally inaccessible $R$ can be dealt with similarly. \ok

\subsection{A separation technique}

When we make computation with the martingale representation property, we often need to extract information about a particular stopping time from an entire stochastic integral. We will need the following technique which separates a stopping time from others in a martingale representation.

Suppose the martingale representation property in $\mathbb{F}$ with a representation process $W$. Then, the $(\mathbb{P},\mathbb{F})$ local martingale $X$ takes all the form $\transp H\centerdot W$ for some $W$-integrable predictable process $H$. We call (any version of) the process $H$ the coefficient of $X$ in its martingale representation with respect to the process $W$. This appellation extends naturally to vector valued local martingales.

\bl\label{single-jump}
Let $R$ be any $\mathbb{F}$ stopping time. Let $\xi\in\mathbf{L}^1(\mathbb{P},\mathcal{F}_{R})$. Let  $H$ denote any coefficient of the $(\mathbb{P},\mathbb{F})$ martingale  $\xi\ind_{[R,\infty)}-(\xi\ind_{[R,\infty)})^{\mathbb{F}\cdot p}$ in its martingale representation with respect to $W$.
\ebe
\item
If $R$ is predictable, the two predictable processes $H$ and $H\ind_{[R]}$ are in the same equivalent class with respect to $W$, whose value is determined by the equation on $\{R<\infty\}$
$$
\transp H_{R} \Delta_{R}W=
\xi-\mathbb{E}[\xi|\mathcal{F}_{R-}].
$$
\item
If $R$ is totally inaccessible, the process $H$ satisfies the equations  on $\{R<\infty\}$
$$
\transp H_{R} \Delta_{R}W=\xi,\ \mbox{ and } \
\transp H_{S} \Delta_{S}W =0 \mbox{ on $\{S\neq R\}$},
$$
for any $\mathbb{F}$ stopping time $S$.
\dbe
\el

\textbf{Proof.} Let us consider only a totally inaccessible stopping time $R$.
In this case, $(\xi\ind_{[R,\infty)})^{\mathbb{F}\cdot p}$ is continuous. Computing the jump at $R$ and at $S$ in the equation$$
\xi\ind_{[R,\infty)}-(\xi\ind_{[R,\infty)})^{\mathbb{F}\cdot p}
=\transp H\centerdot W,
$$
we prove the assertions. \ok

\

\subsection{Representation process reconstituted}\label{decompsection}

As mentioned before, when we make computations under the martingale representation property, we are not restricted to work with the initially given representation process $W$. In choosing suitable representation process, we can render the computations with martingale representation much easier.

\bd
For a multi-dimensional $\mathbb{F}$ local martingale $X$, we say that it has pathwisely orthogonal components, if $[X_i,  X_j]=0$ for $i\neq j$. For a measurable set $\mathtt{A}$, we say that $X$ has pathwisely orthogonal components outside of $\mathtt{A}$, if $\ind_{\mathtt{A}^c}\centerdot[X_i, X_j]=0$ for $i\neq j$.
\ed

Suppose in the rest of this section that the martingale representation property holds in $\mathbb{F}$ with a $d$-dimensional representation process $W$. The following lemma is well-known (cf. \cite{davis}).

\bl
There exists a continuous $d$-dimensional $\mathbb{F}$ local martingale $X'$ which generates the same stable space as that generated by the components of $W^c$, but with pathwisely orthogonal components (some of the components may be identically null). 
\el

Consider the purely discontinuous part $W^d$. We introduce the following notations. Let $(S_n)_{1\leq n<\mathsf{N}^i}$ ($\mathsf{N}^i\leq \infty$) (resp. $(T_n)_{1\leq n<\mathsf{N}^a}$) be a sequence of $(\mathbb{P},\mathbb{F})$ totally inaccessible (resp. strictly positive $(\mathbb{P},\mathbb{F})$ predictable) stopping times such that $[S_n]\cap [S_{n'}]=\emptyset$ for $n\neq n'$ and $\{s\geq 0:\Delta_sW^{di}\neq 0\}\subset\cup_{n\geq 1}[S_n]$ (resp. $[T_n]\cap [T_{n'}]=\emptyset$ for $n\neq n'$ and $\{s\geq 0:\Delta_sW^{da}\neq 0\}\subset\cup_{n\geq 1}[T_n]$).

For $1\leq n'<\mathsf{N}^a$, for $1\leq n<\mathsf{N}^i$, we find and enumerate the partition sets defined in Lemma \ref{partition} for the stopping times $T_{n'}$ or $S_n$: $(A_{n',0},A_{n',1},A_{n',2},\ldots,A_{n',d})$ and $(B_{n,1},B_{n,2},\ldots,B_{n,d})$ (where some $A_{h'}$ and $B_h$ may be empty) such that $$
\mathcal{F}_{T_{n'}}=\mathcal{F}_{T_{n'}-}\vee\sigma(A_{n',0},A_{n',1},A_{n',2},\ldots,A_{n',d})\ \mbox{ and } 
\mathcal{F}_{S_n}=\mathcal{F}_{S_n-}\vee\sigma(B_{n,1},B_{n,2},\ldots,B_{n,d}).
$$
Let $$
p_{n',h'}=\mathbb{P}[A_{n',h'}|\mathcal{F}_{T_{n'}-}]\ \ \mbox{ and }\
q_{n,h}=\mathbb{P}[B_{n,h}|\mathcal{F}_{S_n-}], \ 0\leq h'\leq d, \ 1\leq h\leq d.
$$ 
Let $$
v_{n',h'}\in\mathcal{F}_{T_n-}\ \mbox{ and respectively }\  w_{n,h}\in\mathcal{F}_{S_n-}
$$ 
be the vector value of $\Delta_{T_{n'}} W$ on $A_{n',h'}$ and respectively the vector value of $\Delta_{S_n} W$ on $B_{n,h}$ (cf. Lemma \ref{partition}). We define real processes 
\begin{equation}\label{Xprimeprime}
\dcb
X''_{h'}
=\sum_{1\leq n'<\mathsf{N}^a}\frac{1}{2^n}(\ind_{A_{n',h'}}\ind_{[T_{n'},\infty)}-(\ind_{A_{n',h'}}\ind_{[T_{n'},\infty)})^{\mathbb{F}\cdot p}), \ 0\leq h'\leq d.
\dce
\end{equation}
We define $d$-dimensional vector valued processes
\begin{equation}\label{Xprimeprimeprime}
\dcb
X'''_h
=
\sum_{1\leq n<\mathsf{N}^i}\frac{1}{2^n}(w_{n,h}\ind_{B_{n,h}}\ind_{[S_n,\infty)}-(w_{n,h}\ind_{B_{n,h}}\ind_{[S_n,\infty)})^{\mathbb{F}\cdot p}), \ 1\leq h\leq d,
\dce
\end{equation}
(which is well-defined). Let $X$ be a multi-dimensional local martingale whose components incorporate the processes $X', X'', X'''$.

\bethe\label{pathortho}
The process $X$ has the martingale representation property in $\mathbb{F}$ under $\mathbb{P}$.
\ethe

\textbf{Proof.}
Let $Y$ be a (real) bounded $(\mathbb{P},\mathbb{F})$ martingale orthogonal to the components of $X$. The bracket $[Y,X]$ is a vector valued $(\mathbb{P},\mathbb{F})$ local martingales. By the martingale representation property of $W$ in $(\mathbb{P},\mathbb{F})$, $Y$ takes the form $Y=\transp H\centerdot W$ for some vector valued $\mathbb{F}$ predictable process $H$. The computation of the bracket gives $$
\dcb
[Y,X''_h]&=&\sum_{1\leq n<\mathsf{N}^a}\frac{1}{2^n} [Y,\ \ind_{A_{n,h}}\ind_{[T_n,\infty)}-(\ind_{A_{n,h}}\ind_{[T_n,\infty)})^{\mathbb{F}\cdot p}]\\
&=& \sum_{1\leq n<\mathsf{N}^a}\frac{1}{2^n}\transp H\centerdot [W,\ \ind_{A_{n,h}}\ind_{[T_n,\infty)}-(\ind_{A_{n,h}}\ind_{[T_n,\infty)})^{\mathbb{F}\cdot p}]\\
&=& \sum_{1\leq n<\mathsf{N}^a}\frac{1}{2^n}\transp H_{T_n} \Delta_{T_n}W(\ind_{A_{n,h}}-p_{n,h})\ind_{[T_n,\infty)}\\
&=& \sum_{1\leq n<\mathsf{N}^a}\frac{1}{2^n}\transp H_{T_n} \Delta_{T_n}W\ind_{A_{n,h}}\ind_{[T_n,\infty)}-\sum_{1\leq n<\mathsf{N}^a}\frac{1}{2^n}\transp H_{T_n} \Delta_{T_n}Wp_{n,h}\ind_{[T_n,\infty)}\\

&=& \sum_{1\leq n<\mathsf{N}^a}\frac{1}{2^n}\transp H_{T_n} v_{n,h}\ind_{A_{n,h}}\ind_{[T_n,\infty)}-\sum_{1\leq n<\mathsf{N}^a}\frac{1}{2^n}\transp H_{T_n} \Delta_{T_n}Wp_{n,h}\ind_{[T_n,\infty)}.
\dce
$$ 
It is a $(\mathbb{P},\mathbb{F})$ local martingale. For every $1\leq n<\mathsf{N}^a$, taking the stochastic integral of the predictable process $\ind_{[T_n]}$ with respect to this local martingale, we see that each term
$$
\frac{1}{2^n}\transp H_{T_n} v_{n,h}\ind_{A_{n,h}}\ind_{[T_n,\infty)}-\frac{1}{2^n}\transp H_{T_n} \Delta_{T_n}Wp_{n,h}\ind_{[T_n,\infty)}
$$
is itself a local martingale. Taking the predictable dual projection, we obtain $$
\frac{1}{2^n}\transp H_{T_n}v_{n,h}p_{n,h}\ind_{[T_n,\infty)}\equiv 0 \mbox{ (a  null process)},
$$ 
because $\mathbb{E}[\transp H_{T_n}\Delta_{T_n}W|\mathcal{F}_{T_n}]=0$ and $v_{n,h}\in\mathcal{F}_{T_n-}$. Consequently $\transp H_{T_n}v_{n,h}\ind_{A_{n,h}}=0$ on $\{T_n<\infty\}$. This being true for any $0\leq h\leq d$, we can write $$
\dcb
\Delta_{T_n}Y&=&\transp H_{T_n}\Delta_{T_n}W
=\sum_{h=0}^d\transp H_{T_n}\Delta_{T_n}W\ind_{A_{n,h}}\\
&=&\sum_{h=0}^d\transp H_{T_n}v_{n,h}\ind_{A_{n,h}}
=0,\ \mbox{ on ${\{T_n<\infty\}}$.}
\dce
$$
In the same way,
$$
\dcb
[Y,X'''_h]&=&
\sum_{1\leq n<\mathsf{N}^i}\frac{1}{2^n} [Y,\ w_{n,h}\ind_{B_{n,h}}\ind_{[S_n,\infty)}-(w_{n,h}\ind_{B_{n,h}}\ind_{[S_n,\infty)})^{\mathbb{F}\cdot p}]\\

&=& \sum_{1\leq n<\mathsf{N}^i}\frac{1}{2^n} [Y,\ w_{n,h}\ind_{B_{n,h}}\ind_{[S_n,\infty)}]\
\mbox{ because $(w_{n,h}\ind_{B_{n,h}}\ind_{[S_n,\infty)})^{\mathbb{F}\cdot p}$ is continuous,}\\

&=& \sum_{1\leq n<\mathsf{N}^a}\frac{w_{n,h}}{2^n}\transp H_{S_n} \Delta_{S_n}W\ \ind_{B_{n,h}}\ind_{[S_n,\infty)}\\

&=& \sum_{1\leq n<\mathsf{N}^a}\frac{w_{n,h}}{2^n}\transp H_{S_n} w_{n,h}\ \ind_{B_{n,h}}\ind_{[S_n,\infty)}
\dce
$$
is a local martingale. For $1\leq i<\mathsf{N}^i$, set $J_{i}$ the coefficient of $\ind_{B_{i,h}}\ind_{[S_i,\infty)}-(\ind_{B_{i,h}}\ind_{[R,\infty)})^{\mathbb{F}\cdot p}$ in its martingale representation with respect to $W$. By Lemma \ref{single-jump},
$$
\dcb
&&\transp J_{i}{\centerdot}\left(\sum_{1\leq n<\mathsf{N}^a}\frac{w_{n,h}}{2^n}\transp H_{S_n} w_{n,h}\ \ind_{B_{n,h}}\ind_{[S_n,\infty)}\right)
=
\sum_{1\leq n<\mathsf{N}^a}\frac{\transp J_{i,S_n} w_{n,h}}{2^n}\ \transp H_{S_n} w_{n,h}\ \ind_{B_{n,h}}\ind_{[S_n,\infty)}\\
&=&
\frac{1}{2^i}\ \transp J_{i,S_i} w_{i,h} \transp H_{S_i} w_{i,h}\ \ind_{B_{i,h}}\ind_{[S_i,\infty)}
=
\frac{1}{2^i}\ \transp H_{S_i} w_{i,h}\ \ind_{B_{i,h}}\ind_{[S_i,\infty)}.
\dce
$$
By assumption, it is a $(\mathbb{P},\mathbb{F})$ local martingale. Hence, $\frac{1}{2^i}\ \transp (H_{S_i}w_{n,h}q_{n,h}\ind_{[S_n,\infty)})^{\mathbb{F}\cdot p}$ is a null process. Repeating the reasoning in the preceding paragraphs, we conclude that $\transp H_{S_i} w_{i,h}\ind_{B_{i,h}}=0$ so that $\Delta_{S_i}Y=0$. 

Hence $Y$ is a continuous martingale. As the continuous components $X'$ generate $W^c$, the bracket $\transp H\centerdot\cro{W^c,W^c_k}, 1\leq k\leq d,$ is a local martingale, which must be null. Consequently $ H\centerdot W^c=0$. This proves the theorem, according to \cite[Corollaire(4.12)]{Jacodlivre}. \ok

\

\subsection{The finite predictable constraint condition of representation processes}
\label{rpfpcc}

Under the martingale representation property, the representation processes satisfy the finite predictable constraint. This can result from \cite[Th\'eor\`eme 4.80]{Jacodlivre}. For our account of the martingale representation property to be complete, to follow the logic of the present paper, we give here a brief description of this result, in terms of the processes $X'', X'''$.

\bl\label{pjumpvalue}
There exist a finite number $\mathsf{n}''$ of $d+1$-dimensional $\mathbb{F}$ predictable process $\alpha''_k, 1\leq k\leq \mathsf{n}''$, such that $$
\Delta X'' \ind_{\{\Delta X''\neq \boldsymbol{0}\}}
= 
\sum_{k=1}^{\mathsf{n}''}\alpha''_k\ind_{\{\Delta X''=\alpha''_k \}}.
$$
\el

\textbf{Proof.} 
It is enough to notice that, for every component $$
\dcb
X''_h=\sum_{1\leq n<\mathsf{N}^a}\frac{1}{2^n}(\ind_{A_{n,h}}\ind_{[T_n,\infty)}-(\ind_{A_{n,h}}\ind_{[T_n,\infty)})^{\mathbb{F}\cdot p})
=\sum_{1\leq n<\mathsf{N}^a}\frac{1}{2^n}(\ind_{A_{n,h}}-p_{n,h})\ind_{[T_n,\infty)},
\dce
$$
$0\leq h\leq d$, the jump process $\Delta X''_h$ takes one of the three values $\boldsymbol{0},\zeta_{h,1}, \zeta_{h,2}$, where$$
\dcb
\zeta_{h,1}
=
\sum_{1\leq n<\mathsf{N}^a}\frac{1}{2^n}(1-p_{n,h})\ind_{[T_n]},\
\ \
\zeta_{h,2}
=
\sum_{1\leq n<\mathsf{N}^a}\frac{1}{2^n}(-p_{n,h})\ind_{[T_n]},
\dce
$$
which are $\mathbb{F}$ predictable processes. \ \ok

\bl\label{ijumpvalue}
Suppose that $W$ is locally square integrable. There exist a finite number $\mathsf{n}'''$ of $d\times d$-dimensional $\mathbb{F}$ predictable process $\alpha'''_k, 1\leq k\leq \mathsf{n}'''$, such that$$
\Delta X''' \ind_{\{\Delta X'''\neq \boldsymbol{0}\}}
= 
\sum_{k=1}^{\mathsf{n}'''}\alpha'''_k\ind_{\{\Delta X'''=\alpha'''_k\}}.
$$
\el

\textbf{Proof.} 
Notice that $X'''$ is here locally square integrable. For $1\leq h\leq d$, let $\mu$ be the jump measure of $X'''_h$ with $\mathbb{F}$ compensator $\nu$. We consider, for $1\leq i,j\leq d$, the local martingale $x_ix_j{_*}(\mu-\nu)$ and its representation by $X$. By the pathwise orthogonality, its representation depends only on $X'''_{h}$. Hence, there exists a $\mathbb{F}$ predictable process $H_{j}, 1\leq j\leq d,$ such that $$
\transp H_{j}\centerdot X'''_h
=
x_ix_j{_*}(\mu-\nu).
$$
Computing the jumps, we obtain $$
\sum_{k=1}^d H_{j,k}\Delta X'''_{h,k} = \Delta X'''_{h,i}\Delta X'''_{h,j},\
1\leq j\leq d. 
$$
This means that, if $\Delta X'''_{h}\neq 0$, $\Delta X'''_{h,i}$ is a root of the characteristic polynomial of the matrix of components $(H_{j,k})_{1\leq j,k\leq d}$. Applying \cite[Theorem 2.2]{BS}, there exists $\mathbb{F}$ predictable processes $(\zeta_1,\ldots,\zeta_d)$ such that $$
\Delta X'''_{h,i}\prod_{j=1}^d(\Delta X'''_{h,i}-\zeta_j)=0.
$$
The lemma can now be deduced from this property. \ok

\bethe\label{WY}
Suppose that $W$ has the martingale representation property in $\mathbb{F}$ under $\mathbb{P}$. Then, the process $W$ satisfies the finite $\mathbb{F}$ predictable constraint condition. More precisely, there exist a finite number $\mathsf{n}$ of $d$-dimensional $\mathbb{F}$ predictable process $\alpha_k, 1\leq k\leq \mathsf{n}$, such that$$
\Delta W \ind_{\{\Delta W\neq \boldsymbol{0}\}}
= 
\sum_{k=1}^{\mathsf{n}}\alpha_k\ind_{\{\Delta W=\alpha_k\}}.
$$
\ethe

\textbf{Proof.}
If $W$ is locally square integrable, the theorem is the consequence of the representation of $W$ by $X$ (Theorem \ref{pathortho}) and of Lemma \ref{pjumpvalue} and Lemma \ref{ijumpvalue}. If not, let $T>0$ be a constant. Let $W^*_T=\sup_{s\leq T}|W_s|$. Let $\eta=e^{-W^*_T}$ and let $(\eta_t)_{t\in[0,T]}$ be the associated $(\mathbb{P},\mathbb{F})$ bounded martingale. Let $\overline{\mathbb{P}}=\eta.\mathbb{P}$ and$$
\overline{W} = W - \frac{1}{\eta_-}{\centerdot}\cro{\eta,W}^{\mathbb{P}\cdot\mathbb{F}}.
$$ 
The process $\overline{W}$ is locally square integrable under $\overline{\mathbb{P}}$ on the interval $[0,T]$, and by \cite{JS}, $\overline{W}$ possesses the martingale representation property in $\mathbb{F}$ under $\overline{\mathbb{P}}$. There exist, therefore, a finite number $\mathsf{n}$ (independent of $T$) of $\mathbb{F}$ predictable processes $(\zeta_1,\ldots,\zeta_\mathsf{n})$ such that
$$
\Delta \overline{W} 
\in\{\boldsymbol{0}, \zeta_1,\ldots\ldots,\zeta_\mathsf{n}\}
$$
on $[0,T]$, or equivalently
$$
\Delta W 
\in\{\boldsymbol{0}, \zeta_1+\frac{1}{\eta_-}\Delta\cro{\eta,W}^{\mathbb{P}\cdot\mathbb{F}},\ldots\ldots,\zeta_\mathsf{n}+\frac{1}{\eta_-}\Delta\cro{\eta,W}^{\mathbb{P}\cdot\mathbb{F}}\}
$$
on $[0,T]$. The theorem is deduced from this property. \ok

\bcor\label{ijumpvalue-modif}
Lemma \ref{ijumpvalue} remain available, without the local square integrability of $W$.
\ecor

\

\subsection{Another modification of the representation process}\label{newmodif}

The process $X'''$ is not always locally bounded and has not necessarily pathwisely orthogonal components. With the finite predictable constraint condition, we can modify it to have the boundedness and the pathwise orthogonality.

Consider the process $\alpha'''$ in Lemma \ref{ijumpvalue} (cf. Corollary \ref{ijumpvalue-modif}). Let $\mu$ be the jump measure of $X'''$ with $\mathbb{F}$ compensator $\nu$. Let $e_k, u_k, 1\leq k\leq \mathsf{n}'''$, be the function in Theorem \ref{specif} relative to $\alpha'''$. Let   $$
X^{\circ}_k=u_k{_*}(\mu-\nu),\ 1\leq k\leq \mathsf{n}'''.
$$
The local martingales $X^{\circ}_k$ are mutually pathwisely orthogonal.

\bl\label{modifcirc}
For any $X'''$-integrable $\mathbb{F}$ predictable process $H$, $$
\transp H\centerdot X'''
=
\transp H\frac{\alpha'''}{\mathsf{e}(\alpha''')}\ind_{\{\mathsf{e}(\alpha''')\neq 0\}} \centerdot X^{\circ}.
$$
\el

\textbf{Proof.}
We have, by \cite[Theorem 11.23 and 11.24]{HWY}, $$
\transp H\centerdot X'''
=
\transp H\centerdot (x{_*}(\mu-\nu))
=
\transp Hx{_*}(\mu-\nu).
$$
Applying Theorem \ref{specif}, we see that $\transp H\frac{\alpha'''}{\mathsf{e}(\alpha''')}\ind_{\{\mathsf{e}(\alpha''')\neq 0\}}$ is ${X^{\circ}}$-integrable and
$$
\transp H\centerdot X'''
=
\transp H\frac{\alpha'''}{\mathsf{e}(\alpha''')}\ind_{\{\mathsf{e}(\alpha''')\neq 0\}} \centerdot {X^{\circ}}.\ \ok
$$

\bcor\label{XXXo}
The process $(X',X'',{X^{\circ}})$ possesses the martingale representation property in $\mathbb{F}$ under $\mathbb{P}$.
\ecor

\brem\label{allorthog}
We note that $(X',X'',{X^{\circ}})$ is a locally bounded process. The three processes $X',X'',{X^{\circ}}$ are mutually pathwisely orthogonal. The components of the processes $X',{X^{\circ}}$ are pathwisely orthogonal. The path of $X''$ is of finite variation. Let $H$ be a $(X',X'',X^\circ)$-integrable predictable process. The process $H$ is naturally cut into three parts $(H',H'',H''')$ corresponding to $(X',X'',X^\circ)$. By the pathwise orthogonality, $H'_h$ is $X'_h$-integrable for $1\leq h\leq d$, $H''$ is $X''$-integrable, $H'''_h$ is $X^\circ_h$-integrable for $1\leq h\leq \mathsf{n}'''$.
\erem

\bethe\label{boundedW}
If the martingale representation property holds in $\mathbb{F}$ under $\mathbb{P}$, there exists always a locally bounded representation process, which has pathwisely orthogonal components outside of a predictable thin set.
\ethe

\

\section{Fully viable market expansion and the drift multiplier assumption}\label{dma}

\subsection{The setting}

Let $\mathbb{G}=(\mathcal{G}_t)_{t\geq 0}$ be a second filtration of sub-$\sigma$-algebras of $\mathcal{A}$, containing $\mathbb{F}$, i.e., $\mathcal{G}_t\supset\mathcal{F}_t$ for $t\geq 0$. We call $\mathbb{G}$ an expansion (or an enlargement) of the filtration $\mathbb{F}$.

\subsubsection{Local martingale deflator}

We recall the notion of deflator.

\bd\label{scdef}
Let $T$ be a $\mathbb{G}$ stopping time. We call a strictly positive $\mathbb{G}$ adapted real process $Y$ with $Y_0=1$, a local martingale deflator on the time horizon $[0,T]$ for a (multi-dimensional) $(\mathbb{P},\mathbb{G})$ special semimartingale $X$, if the processes $Y$ and $Y X$ are $(\mathbb{P},\mathbb{G})$ local martingales on $[0,T]$. 
\ed

This is a notion of no-arbitrage condition. Actually, the existence of local martingale deflators is equivalent to the no-arbitrage conditions \texttt{NUPBR} and \texttt{NA}$_{1}$ (cf. \cite{kabanov, KC2010, song-takaoka, takaoka}). We know that, when the no-arbitrage condition \texttt{NUPBR} is satisfied, the market is viable, and vice versa.

\

\subsubsection{Full viability}

We consider in this section the following assumption. Let $T$ be a $\mathbb{G}$ stopping time.

\bassumption\label{fullvaibility}(\textbf{Full viability} on $[0,T]$)
The expansion from $\mathbb{F}$ to $\mathbb{G}$ is fully viable on $[0,T]$. This means that, for any strictly positive $\mathbb{F}$ local martingale $X$, $X$ has the no-arbitrage property of the first kind in $\mathbb{G}$ on $[0,T]$, i.e. $X$ has a (local martingale) deflator in $\mathbb{G}$ on $[0,T]$. 
\eassumption

We refer to \cite{kabanov, KC2010, K2012}, also to Definition \ref{scdef} below, for the notion of no-arbitrage of the first kind and the notion of deflator. We can check that, if $\mathbb{G}$ is fully viable, the hypothesis$(H')$ (cf. \cite{Jacod, J}) from $\mathbb{F}$ to $\mathbb{G}$ is satisfied.

\bassumption
\textit{
(\noindent\textbf{Hypothesis$(H')$} on the time horizon $[0,T]$) 
}
Every $(\mathbb{P},\mathbb{F})$ local martingale is a $(\mathbb{P},\mathbb{G})$ semimartingale on $[0,T]$.
\eassumption

Whenever Hypothesis$(H')$ holds, the associated drift operator can be defined.

\bl\label{linearG}
Suppose hypothesis$(H')$ on $[0,T]$. Then there exists a linear map $\Gamma$ from the space of all $(\mathbb{P},\mathbb{F})$ local martingales into the space of c\`adl\`ag $\mathbb{G}$-predictable processes on $[0,T]$, with finite variation and null at the origin, such that, for any $(\mathbb{P},\mathbb{F})$ local martingale $X$, $\widetilde{X}:=X-\Gamma(X)$ is a $(\mathbb{P},\mathbb{G})$ local martingale on $[0,T]$. Moreover, if $X$ is an $\mathbb{F}$ local martingale and $H$ is an $\mathbb{F}$ predictable $X$-integrable process, then $H$ is $\Gamma(X)$-integrable and $\Gamma(H\centerdot X)=H\centerdot \Gamma(X)$ on $[0,T]$. The operator $\Gamma$ will be called the drift operator.
\el

\textbf{Proof.} Note that, under hypothesis$(H')$, for any $\mathbb{F}$ local martingale $X$, $X$ is a special $\mathbb{G}$ semimartingale on $[0,T]$ (cf. \cite[Definition 8.4 and Theorem 8.6]{HWY}) so that the drift operator is well-defined. The linearity of $\Gamma$ is the consequence of the uniqueness of special semimartingale decomposition (cf. \cite[Theorem 8.5]{HWY}). The property of $\Gamma(H\centerdot X)$ is the consequence of \cite[Lemma 2.2]{JS}. \ok

\

\subsubsection{Drift multiplier assumption}

But in this paper we are actually interested in the following extra assumption.

\bassumption\label{assump1} (\textbf{Drift multiplier assumption} on $[0,T]$) Let $T$ be a $\mathbb{G}$ stopping time.
\
\ebe 
\item
The Hypothesis$(H')$ is satisfied on the time horizon $[0,T]$ with a drift operator $\Gamma$.

\item
There exist $N=(N_1,\ldots,N_\mathsf{n})$ an $\mathsf{n}$-dimensional $(\mathbb{P},\mathbb{F})$ local martingale, and ${\varphi}$ an $\mathsf{n}$ dimensional $\mathbb{G}$ predictable process such that, for any $(\mathbb{P},\mathbb{F})$ local martingale $X$, $[N,X]^{\mathbb{F}\cdot p}$ exists, ${\varphi}$ is $[N,X]^{\mathbb{F}\cdot p}$-integrable, and $$
\Gamma(X)=\transp{\varphi}\centerdot [N,X]^{\mathbb{F}-p}
$$
on the time horizon $[0,T]$.

\dbe
\eassumption

This section is devoted the proof of the following theorem.

\bethe\label{driftsatisfied}
Suppose the martingale representation property in $\mathbb{F}$. Suppose the full viability on $[0,T]$. Then, $\Gamma$ satisfies the drift multiplier assumption on $[0,T]$.  
\ethe

\textbf{Proof.}
It is the consequence of Lemma \ref{driftprime1}, Lemma \ref{driftprimeprime}, Lemma \ref{driftprime3} in the next subsections,  with help of Lemma \ref{linearG} and Remark \ref{allorthog}. In fact, under  the martingale representation property in $\mathbb{F}$, we have a representation process $(X',X'',X^\circ)$ defined in Corollary \ref{XXXo}. For any $\mathbb{F}$ local martingale $X$, we write $X$ in its representation $H'\centerdot X'+H''\centerdot X''+H^\circ\centerdot X^\circ$ with respect to the representation process $(X',X'',X^\circ)$ (cf. Remark \ref{allorthog}). Let $\mathtt{A}$ be any $\mathbb{F}$ predictable set such that $H'\ind_{\mathtt{A}}, H''\ind_{\mathtt{A}}, H^\circ\ind_{\mathtt{A}}$ are bounded. By Lemma \ref{driftprime1}, Lemma \ref{driftprimeprime}, Lemma \ref{driftprime3}, Lemma \ref{linearG} and Remark \ref{allorthog}, we compute, on $[0,T]$,$$
\dcb
&&\ind_{\mathtt{A}}\centerdot\Gamma(X)
=
\Gamma(\ind_{\mathtt{A}}\centerdot X)
=
\Gamma(\ind_{\mathtt{A}}H'\centerdot X')
+
\Gamma(\ind_{\mathtt{A}}H''\centerdot X'')
+
\Gamma(\ind_{\mathtt{A}}H^\circ\centerdot X^\circ)\\

&=&
\sum_{h}\Gamma(\ind_{\mathtt{A}}H'_{h}\centerdot X'_{h})
+
\Gamma(\ind_{\mathtt{A}}H''\centerdot X'')
+
\sum_{h}\Gamma(\ind_{\mathtt{A}}H^\circ_{h}\centerdot X^\circ_{h})\\

&=&
\sum_{h}\ind_{\mathtt{A}}H'_{h}\centerdot \Gamma(X'_{h})
+
\ind_{\mathtt{A}}H''\centerdot \Gamma(X'')
+
\sum_{h}\ind_{\mathtt{A}}H^\circ_{h}\centerdot \Gamma(X^\circ_{h})\\

&=&
\sum_{h}\ind_{\mathtt{A}}H'_{h}  G'_{h}\centerdot[X'_h,X'_h]^{\mathbb{F}\cdot p}
+
\ind_{\mathtt{A}}H''  \transp\varphi''\centerdot[N'',X''_h]^{\mathbb{F}\cdot p}
+
\sum_{h}\ind_{\mathtt{A}}H^\circ_{h}  G^\circ_{h}\centerdot [X^\circ_h, X^\circ_h]^{\mathbb{F}\cdot p}\\

&=&
\sum_{h}  \ind_{\mathtt{A}}G'_{h}\centerdot[X'_h, H'_{h}\centerdot X'_h]^{\mathbb{F}\cdot p}
+
\ind_{\mathtt{A}}\transp\varphi''\centerdot[N'', H''_{h}\centerdot X''_h]^{\mathbb{F}\cdot p}
+
\sum_{h} \ind_{\mathtt{A}}G^\circ_{h}\centerdot [X^\circ_h, H^\circ_{h}\centerdot X^\circ_h]^{\mathbb{F}\cdot p}\\

&=&
\sum_{h} \ind_{\mathtt{A}} G'_{h}\centerdot[X'_h, X]^{\mathbb{F}\cdot p}
+
\ind_{\mathtt{A}}\transp\varphi''\centerdot[N'', X]^{\mathbb{F}\cdot p}
+
\sum_{h} \ind_{\mathtt{A}}G^\circ_{h}\centerdot [X^\circ_h, X]^{\mathbb{F}\cdot p}.
\dce
$$
This proves first of all the necessary integrability conditions and then the theorem. \
\ok

\

\subsection{General consequences of the full viability}\label{gcfv}

We begin with an immediate consequence of the full viability on the drift operator.

\bl\label{fbd}
Let $T$ be a $\mathbb{G}$ stopping time. Suppose that the expansion is fully viable on $[0,T]$. For any $\mathbb{F}$ locally bounded local martingale $X$, there exists a strictly positive $\mathbb{G}$ local martingale $Y$ such that $$
\Gamma(X)
=
-\frac{1}{Y_-}\centerdot[Y,X]^{\mathbb{G}\cdot p}\ \mbox{ on $[0,T]$}.
$$
\el

\textbf{Proof.}
For any $\mathbb{F}$ stopping time $T'$ such that $X^{T'}$ is bounded, for some $a>0$, $a|\Delta X|<1$ on $[0,T']$. Let $S=\mathcal{E}(aX)$ which is strictly positive. By the full viability, there exists a strictly positive $\mathbb{G}$ local martingale $Y$ such that $YS$ is a $\mathbb{G}$ local martingale on $[0,T'\wedge T]$. The lemma is the consequence of the integration by parts formula$$
\dcb
YS=Y_0S_0+S_-\centerdot Y + Y_-\centerdot S + [Y,S]\ \mbox{ or equivalently}\\
aY_-S_-\centerdot X+aS_{-}\centerdot[Y,X]
=
YS-Y_0S_0-S_-\centerdot Y \ \mbox{ on $[0,T'\wedge T]$.}\ \ok
\dce
$$ 

\bl\label{interm}
Suppose the full viability of the expansion on $[0,T]$. For any $\mathbb{F}$ locally bounded $\mathbb{F}$ optional process $A$ with finite variation, there exists a strictly positive $\mathbb{G}$ local martingale $Y$ such that $$
(Y\centerdot A)^{\mathbb{G}\cdot p}
=
Y_-\centerdot A^{\mathbb{F}\cdot p}\ \mbox{ on $[0,T]$.}
$$
Consequently $A^{\mathbb{G}\cdot p}$ is absolutely continuous with respect to $A^{\mathbb{F}\cdot p}$ on $[0,T]$.
\el

\textbf{Proof.}
For any $\mathbb{F}$ stopping time $T'>0$ such that $A^{T'}$ is bounded, for some $a>0$, $S=\mathcal{E}(a(A-A^{\mathbb{F}\cdot p}))$ is strictly positive on $[0,T']$. There exists a strictly positive $\mathbb{G}$ local martingale $Y$ such that $YS$ is a $\mathbb{G}$ local martingale on $[0,T'\wedge T]$. Write the integration by parts formula$$
\dcb
YS=Y_0S_0+S_-\centerdot Y + Y\centerdot S\ \mbox{ or equivalently }\\
aYS_{-}\centerdot A-aYS_-\centerdot A^{\mathbb{F}\cdot p}
=
YS-Y_0S_0-S_-\centerdot Y
\dce
$$ 
on $[0,T'\wedge T]$. Consequently,$$
(Y\centerdot A)^{\mathbb{G}\cdot p}
=
(Y\centerdot A^{\mathbb{F}\cdot p})^{\mathbb{G}\cdot p}
=
{^{\mathbb{G}\cdot p}}\!(Y)\centerdot A^{\mathbb{F}\cdot p}
=
Y_-\centerdot A^{\mathbb{F}\cdot p}
$$
on $[0,T'\wedge T]$. \ \ok

We now apply Lemma \ref{interm} to view the $\mathbb{F}$ totally inaccessible stopping times. 

\bcor\label{totinacc}
Suppose the full viability of the expansion on $[0,T]$. For any $\mathbb{F}$ totally inaccessible stopping time $S$, there exists a strictly positive $\mathbb{G}$ local martingale $Y$ such that $$
(Y_S\ind_{[S,\infty)})^{\mathbb{G}\cdot p}
=
Y_-\centerdot (\ind_{[S,\infty)})^{\mathbb{F}\cdot p}
$$
on $[0,T]$. Consequently $(\ind_{[S,\infty)})^{\mathbb{G}\cdot p}$ is absolutely continuous with respect to $(\ind_{[S,\infty)})^{\mathbb{F}\cdot p}$ on $[0,T]$, and $S_{\{S\leq T\}}$ is $\mathbb{G}$ totally inaccessible. 
\ecor

\

\subsection{Drift of $X'$}

We suppose for the rest of this paper the martingale representation property in $\mathbb{F}$ with a $d$-dimensional representation process and the full viability of the expansion on $[0,T]$. We use the reconstituted representation process $(X',X'',X^\circ)$ in Corollary \ref{XXXo}. We will compute successively the drifts $\Gamma(X'), \Gamma(X'')$ and $\Gamma(X^\circ)$. 

We begin with the computation of $\Gamma(X')$ which is simple.

\bl\label{driftprime1}
For $1\leq h\leq d$, there exists a $\mathbb{G}$ predictable process $G'_{h}$ such that
$$
\dcb
&&\Gamma(X'_h)
=
G'_{h}\centerdot[X'_h,X'_h]^{\mathbb{F}\cdot p} \ \mbox{ on $[0,T]$}.
\dce
$$
\el

\textbf{Proof.}
Let $Y$ be defined in Lemma \ref{fbd} for $X'_h$. By the path continuity, by \cite{anselstricker}, there exists a $\mathbb{G}$ predictable process $H$ such that $$
[Y,X'_h]^{\mathbb{G}\cdot p} = [H\centerdot \widetilde{X}'_h, \widetilde{X}'_h]
= H\centerdot [\widetilde{X}'_h, \widetilde{X}'_h]
= H\centerdot [X'_h, X'_h]
= H\centerdot [X'_h, X'_h]^{\mathbb{F}\cdot p}
$$
on $[0,T]$. Applying Lemma \ref{fbd}, we prove the lemma. \ok

\

\subsection{Drift of $X''$}

The case of $\Gamma(X'')$ is much more involved. Recall the $\mathbb{F}$ predictable stopping time $T_n$ and the partition sets $A_{n,h}$ defined in subsection \ref{decompsection} for $1\leq n<\mathsf{N}^a, 0\leq h\leq d$.
Let $\mu$ denote the jump measure of $X''$. Let $\nu$ (resp. $\overline{\nu}$) be the $\mathbb{F}$ (resp. $\mathbb{G}$) compensator of $\mu$. We consider the stochastic $_*$-integral in $\mathbb{F}$ with respect to $(\mu-\nu)$, but also in $\mathbb{G}$ with respect to $(\mu-\overline{\nu})$.

\bl\label{Hmnp}
We have $\mathscr{G}^{}(\mathbb{F},\mu)\ind_{[0,T]}\subset\mathscr{G}^{}(\mathbb{G},\mu)$. For $g\in\mathscr{G}(\mathbb{F},\mu)$, on $[0,T]$, $$
g{_*}(\mu-\overline{\nu})
=
g{_*}(\mu-{\nu})-\Gamma(g{_*}(\mu-{\nu})),
$$
and
$$
\Gamma(g{_*}(\mu-{\nu}))
=
\sum_{n=1}^{\mathsf{N}^a-}\left(\ \mathbb{E}[g(T_n,\Delta_{T_n}X'')\ind_{\{\Delta_{T_n}X''\neq 0\}}|\mathcal{G}_{T_n-}]
-\mathbb{E}[g(T_n,\Delta_{T_n}X'')\ind_{\{\Delta_{T_n}X''\neq 0\}}|\mathcal{F}_{T_n-}]\ \right)\ind_{[T_n\infty)}.
$$
In particular, $x{_*}(\mu-\overline{\nu})=\widetilde{X}''$ on $[0,T]$.
\el 

\textbf{Proof.}
For $1\leq n<\mathsf{N}^a$, the process $$
(g(T_n,\Delta_{T_n}X'')\ind_{\{\Delta_{T_n}X''\neq 0\}}-\mathbb{E}[g(T_n,\Delta_{T_n}X'')\ind_{\{\Delta_{T_n}X''\neq 0\}}	|\mathcal{F}_{T_n-}])\ind_{[T_n\infty)}
=
\ind_{[T_n]}\centerdot(g{_*}(\mu-\nu))
$$
is a $\mathbb{F}$ local martingale. The martingale part in $\mathbb{G}$ of this process is given by
$$
\dcb
&&
(g(T_n,\Delta_{T_n}X'')\ind_{\{\Delta_{T_n}X''\neq 0\}}-\mathbb{E}[g(T_n,\Delta_{T_n}X'')\ind_{\{\Delta_{T_n}X''\neq 0\}}|\mathcal{F}_{T_n-}])\ind_{[T_n\infty)}\\
&&
-\mathbb{E}[(g(T_n,\Delta_{T_n}X'')\ind_{\{\Delta_{T_n}X''\neq 0\}}-\mathbb{E}[g(T_n,\Delta_{T_n}X'')\ind_{\{\Delta_{T_n}X''\neq 0\}}|\mathcal{F}_{T_n-}])|\mathcal{G}_{T_n-}]\ind_{[T_n\infty)}\\
&=&
(g(T_n,\Delta_{T_n}X'')\ind_{\{\Delta_{T_n}X''\neq 0\}}-\mathbb{E}[g(T_n,\Delta_{T_n}X'')\ind_{\{\Delta_{T_n}X''\neq 0\}}|\mathcal{G}_{T_n-}])\ind_{[T_n\infty)}.
\dce
$$
(In particular, this shows that $\mathbb{E}[g(T_n,\Delta_{T_n}X'')\ind_{\{\Delta_{T_n}X''\neq 0\}}|\mathcal{G}_{T_n-}]$ is well-defined.)
By Lemma \ref{linearG}, this implies also$$
\dcb
&&\ind_{[T_n]}\centerdot\Gamma(g{_*}(\mu-\nu))
=
\Gamma(\ind_{[T_n]}\centerdot(g{_*}(\mu-\nu)))\\
&=&
(\mathbb{E}[g(T_n,\Delta_{T_n}X'')\ind_{\{\Delta_{T_n}X''\neq 0\}}|\mathcal{G}_{T_n-}]
-\mathbb{E}[g(T_n,\Delta_{T_n}X'')\ind_{\{\Delta_{T_n}X''\neq 0\}}|\mathcal{F}_{T_n-}])\ind_{[T_n\infty)}
\dce
$$
on $[0,T]$. Because $\Gamma(g{_*}(\mu-\nu))$ is $\mathbb{G}$ predictable with finite variation on $[0,T]$, the series$$
\sum_{n=1}^{\mathsf{N}^a-}\left|\ \mathbb{E}[g(T_n,\Delta_{T_n}X'')\ind_{\{\Delta_{T_n}X''\neq 0\}}|\mathcal{G}_{T_n-}]
-\mathbb{E}[g(T_n,\Delta_{T_n}X'')\ind_{\{\Delta_{T_n}X''\neq 0\}}|\mathcal{F}_{T_n-}]\ \right|\ind_{[T_n\infty)}
$$
is a $\mathbb{G}$ locally integrable predictable process on $[0,T]$. This local integrability, together with the relation $$
\dcb
g(T_n,\Delta_{T_n}X'')\ind_{\{\Delta_{T_n}X''\neq 0\}}-
\mathbb{E}[g(T_n,\Delta_{T_n}X'')\ind_{\{\Delta_{T_n}X''\neq 0\}}|\mathcal{G}_{T_n-}]\\
=
g(T_n,\Delta_{T_n}X'')\ind_{\{\Delta_{T_n}X''\neq 0\}}-
\mathbb{E}[g(T_n,\Delta_{T_n}X'')\ind_{\{\Delta_{T_n}X''\neq 0\}}|\mathcal{F}_{T_n-}]\\
\hspace{3cm} \ +
\mathbb{E}[g(T_n,\Delta_{T_n}X'')\ind_{\{\Delta_{T_n}X''\neq 0\}}|\mathcal{F}_{T_n-}]
-\mathbb{E}[g(T_n,\Delta_{T_n}X'')\ind_{\{\Delta_{T_n}X''\neq 0\}}|\mathcal{G}_{T_n-}],
\dce
$$
implies $g\ind_{[0,T]}\in\mathscr{G}(\mathbb{G},\mu)$, and also$$
\dcb
&&\Gamma(g{_*}(\mu-{\nu}))
=
\Gamma(\ind_{\cup_{1\leq n<\mathsf{N}^a}[T_n]}\centerdot(g{_*}(\mu-{\nu})))
=
\ind_{\cup_{1\leq n<\mathsf{N}^a}[T_n]}\centerdot\Gamma(g{_*}(\mu-{\nu}))\\
&=&
\sum_{n=1}^{\mathsf{N}^a-}\left(\ \mathbb{E}[g(T_n,\Delta_{T_n}X'')\ind_{\{\Delta_{T_n}X''\neq 0\}}|\mathcal{G}_{T_n-}]
-\mathbb{E}[g(T_n,\Delta_{T_n}X'')\ind_{\{\Delta_{T_n}X''\neq 0\}}|\mathcal{F}_{T_n-}]\ \right)\ind_{[T_n\infty)}
\dce
$$
on $[0,T]$. We can now check that $g{_*}(\mu-\overline{\nu})$ and $g{_*}(\mu-{\nu})-\Gamma(g{_*}(\mu-{\nu}))$ have the same jumps on $[0,T]$. By \cite[Theorem 7.23]{HWY}, they are the same $\mathbb{G}$ local martingales on $[0,T]$. \ok

\brem
Note that, for $0\leq h\leq d$, $X''_h$ is a bounded process with finite variation. $X''_h$ is always a $\mathbb{G}$ special semimartingale whatever hypothesis$(H')$ is valid or not. Denote always by $\widetilde{X}''$ the $\mathbb{G}$ martingale part of $X''$.
\erem

\bl\label{driftprimeprime}
There exist a $d$-dimensional $\mathbb{F}$ local martingale $N''$ of the form $N''=H\centerdot X''$, and a $d$-dimensional $\mathbb{G}$ predictable process $\varphi''$ such that, for every $0\leq h\leq d$, $\transp\varphi''\centerdot[N'',X''_h]^{\mathbb{F}\cdot p}$ exists and$$
\widetilde{X}''_h
=
X''_h-\transp\varphi''\centerdot[N'',X''_h]^{\mathbb{F}\cdot p}
$$
is a $\mathbb{G}$ local martingale. In particular, in case of the full viability on $[0,T]$, $\Gamma(X''_h)=\transp\varphi''\centerdot[N'',X''_h]^{\mathbb{F}\cdot p}$ on $[0,T]$.
\el

\textbf{Proof.}
In this proof, we will simply write $N,\varphi$ instead of $N'',\varphi''$. With the computations in the proof of Lemma \ref{Hmnp}, we know that the $\mathbb{G}$ drift part of $X''_h$ is given by
$$
\sum_{n=1}^{\mathsf{N}^a-} \mathbb{E}[\Delta_{T_n}X''_h|\mathcal{G}_{T_n-}]\ind_{[T_n\infty)}
=
\sum_{n=1}^{\mathsf{N}^a-}\frac{1}{2^n}\left(\ \mathbb{E}[\ind_{A_{n,h}}|\mathcal{G}_{T_n-}]-p_{n,h}\ \right)\ind_{[T_n\infty)}
=
\sum_{n=1}^{\mathsf{N}^a-}\frac{1}{2^n}\left(\ \overline{p}_{n,h}-p_{n,h}\ \right)\ind_{[T_n\infty)}.
$$
with $p_{n,h}=\mathbb{E}[\ind_{A_{n,h}}|\mathcal{F}_{T_n-}]$ and $\overline{p}_{n,h}=\mathbb{E}[\ind_{A_{n,h}}|\mathcal{G}_{T_n-}]$. We look for a $d$-dimensional $\mathbb{F}$ local martingale $N$ and a $d$-dimensional $\mathbb{G}$ predictable process $\varphi$ such that 
$$
\frac{1}{2^n}(\overline{p}_{n,h}-{p}_{n,h})
=
\transp\varphi_{T_n} n_{n,h} p_{n,h},
$$
on $\{T_n<\infty\}$, where $n_{n,h}\in\mathcal{F}_{T_n-}$ is the value of $\Delta_{T_n}N$ on $A_{n,h}$, or equivalently,
\begin{equation}\label{ppfin}
\frac{1}{2^n}(\frac{\overline{p}_{n,h}}{{p}_{n,h}}-1)
=
\transp\varphi_{T_n} n_{n,h},\ 0\leq h\leq d,
\end{equation} 
(with the convention that $\frac{0}{0}-1=0$). Consider the $(1+d)$-dimensional vector $p_n=(p_{n,h})_{0\leq h\leq d}$. By Gram-Schmidt process, we obtain a $\mathcal{F}_{T_n-}$ measurable orthonormal basis $(\epsilon_{n,0},\epsilon_{n,1},\epsilon_{n,2},\ldots,\epsilon_{n,d})$ in $\mathbb{R}\times\mathbb{R}^d$ such that $\epsilon_{n,j}$ is orthogonal to $p_n$ for all $1\leq j\leq d$. Note that
$$
\frac{1}{2^n}\sum_{h=0}^d(\frac{\overline{p}_{n,h}}{{p}_{n,h}}-1)p_{n,h}
=
\frac{1}{2^n}\sum_{h=0}^d({\overline{p}_{n,h}}-p_{n,h})=0.
$$
This implies that the vector $\frac{1}{2^n}(\frac{\overline{\mathsf{p}}}{{\mathsf{p}}}-1)$ of the components $\frac{1}{2^n}(\frac{\overline{p}_{n,h}}{{p}_{n,h}}-1)$ is orthogonal to $p_n$ so that it is a linear combination of the $\epsilon_{n,j}, 1\leq j\leq d$:
$$
\frac{1}{2^n}(\frac{\overline{\mathsf{p}}}{{\mathsf{p}}}-1) 
= \varsigma_{n,1}\epsilon_{n,1}+\varsigma_{n,2}\epsilon_{n,2}+\ldots+\varsigma_{n,d}\epsilon_{n,d},
$$
where $\varsigma_{n,h}$ are the scalar product of $\frac{1}{2^n}(\frac{\overline{\mathsf{p}}}{{\mathsf{p}}}-1)$ with $\epsilon_{n,h}$ so that $\mathcal{G}_{T_n-}$ measurable. Let $\transp\epsilon_n$ denote the $d\times(1+d)$-matrix whose lines are the vectors $\transp\epsilon_{n,k}, 1\leq k\leq d$. Let $\varphi_n$ denote the vector in $\mathbb{R}^d$ of components $2^n\varsigma_{n,k}, 1\leq k\leq d$. Let $n_{n,h}$ denote the $h$th-column ($0\leq h\leq d$) of the matrix $\frac{1}{2^n}\transp \epsilon_n$, which is a vector in $\mathbb{R}^d$. Then, the above identity becomes$$
\dcb
\frac{1}{2^n}\transp(\frac{\overline{\mathsf{p}}}{{\mathsf{p}}}-1) 
= 
\sum_{k=1}^d\varsigma_{n,k}\transp\epsilon_{n,k}
=
\sum_{k=1}^d 2^n\varsigma_{n,k}\ \frac{1}{2^n}\transp\epsilon_{n,k},\\
\mbox{ or }\\
\frac{1}{2^n}(\frac{\overline{p}_{n,h}}{{p}_{n,h}}-1)	
=
\transp\varphi_{n} n_{n,h},\ 0\leq h\leq d.
\dce
$$
The equation (\ref{ppfin}) is solved. We define $d$ number of $\mathbb{F}$ local martingales.$$
N_j=\transp(\sum_{n=1}^{\mathsf{N}^a-}\epsilon_{n,j}\ind_{[T_n]})\centerdot X'',\ 1\leq j\leq d.	
$$
Let $\mathsf{a}_n$ denote the vector of the components $\ind_{A_{n,h}}, 0\leq h\leq d$. We compute the jumps at $T_n<\infty$.$$
\Delta_{T_n}N_j
=\transp \epsilon_{n,j} \Delta_{T_n}X''
=\frac{1}{2^n}\transp \epsilon_{n,j} (\mathsf{a}_n-p_n)
=\frac{1}{2^n}\transp \epsilon_{n,j} \mathsf{a}_n, \ 1\leq j\leq d.
$$
Hence, if  $\varphi=\sum_{n=1}^{\mathsf{N}^a-}2^n\varphi_{n}\ind_{[T_n]}$, for $0\leq h\leq d$,  $$
\dcb
&&\transp\varphi\centerdot[N,X''_h]^{\mathbb{F}\cdot p}
=
\sum_{n=1}^{\mathsf{N}^a-}2^n\transp\varphi_{n}\mathbb{E}[\Delta_{T_n}N\Delta_{T_n}X''_h|\mathcal{F}_{T_n-}]\ind_{[T_n,\infty)}\\

&=&
\sum_{n=1}^{\mathsf{N}^a-}2^n\transp\varphi_{n}\mathbb{E}[\frac{1}{2^n}\transp\epsilon_{n}\mathsf{a}_n\Delta_{T_n}X''_h|\mathcal{F}_{T_n-}]\ind_{[T_n,\infty)}\\

&=&
\sum_{n=1}^{\mathsf{N}^a-}2^n\sum_{j=0}^d\transp\varphi_{n}\mathbb{E}[\frac{1}{2^n}\transp\epsilon_{n}\mathsf{a}_n\Delta_{T_n}X''_h\ind_{A_{n,j}}|\mathcal{F}_{T_n-}]\ind_{[T_n,\infty)}\\

&=&
\sum_{n=1}^{\mathsf{N}^a-}2^n\sum_{j=0}^d\transp\varphi_{n}n_{n,j}x_{n,j}\mathbb{E}[\ind_{A_{n,j}}|\mathcal{F}_{T_n-}]\ind_{[T_n,\infty)}\\

&&
\mbox{where $x_{n,j}=\frac{1}{2^n}(\delta_{h,j}-p_{n,h})$ is the value of $\Delta_{T_n}X''_h$ on $A_{n,j}$,}
\\
&=&
\sum_{n=1}^{\mathsf{N}^a-}2^n\sum_{j=0}^d\transp\varphi_{n}n_{n,j}x_{n,j}p_{n,j}\ind_{[T_n,\infty)}\\

&=&
\sum_{n=1}^{\mathsf{N}^a-}2^n\sum_{j=0}^d \frac{1}{2^n}(\frac{\overline{p}_{n,j}}{{p}_{n,j}}-1)x_{n,j}p_{n,j}\ind_{[T_n,\infty)}\\

&=&
\sum_{n=1}^{\mathsf{N}^a-}\sum_{j=0}^d \frac{1}{2^n}(\overline{p}_{n,j}-p_{n,j})(\delta_{h,j}-p_{n,h})\ind_{[T_n,\infty)}\\

&=&
\sum_{n=1}^{\mathsf{N}^a-} \frac{1}{2^n}(\overline{p}_{n,h}-p_{n,h})\ind_{[T_n,\infty)}\\

&=&
\mbox{ $\mathbb{G}$ drift part of $X''_h$}.
\dce
$$
(Modifying a little the above computation, we can prove that the stochastic integral $\transp\varphi\centerdot[N,X''_h]^{\mathbb{F}\cdot p}$ exists.) \ok

The following lemma will not be used in this paper, but useful in \cite{song-drift}.

\bl\label{covariance-matrix}
For $1\leq n<\mathsf{N}^a$, let $\mathtt{I}_n=\{0\leq h\leq d: p_{n,h}>0\}$. The kernel of the matrix $\Delta_{T_n}[{X}'',\transp {X}'']^{\mathbb{F}\cdot p}
$
on $\{T_n<\infty\}$ is $$
\{a\in\mathbb{R}\times\mathbb{R}^d: \ \mbox{$a_h$ is constant on $h\in\mathtt{I}_n$}\}.
$$ 
There exists a $\mathbb{G}$ predictable matrix valued process $G$ such that $[\widetilde{X}'',\transp \widetilde{X}'']^{\mathbb{G}\cdot p}=G\centerdot [{X}'',\transp X'']^{\mathbb{F}\cdot p}$.
\el

\textbf{Proof.}
Fix $1\leq n<\mathsf{N}^a$. $\mathtt{I}_n$ is an $\mathcal{F}_{T_n-}$ measurable random variable. For an example, suppose $\mathtt{I}_n=\{0,\ldots,k\}$. Let $\mathsf{a}$ denote the vector of the $\ind_{A_{n,h}}, 0\leq h\leq k$, and $p$ denote the vector of the $p_{n,h}, 0\leq h\leq k$. We write
$$
\dcb
&&\left(\mathbb{E}[\Delta_{T_n}{X}''_i\Delta_{T_n}{X}''_j|\mathcal{F}_{T_n-}]\right)_{1\leq i,j\leq k}
=
\frac{1}{4^n}\mathbb{E}[(\mathsf{a}-{p})\transp(\mathsf{a}-{p})|\mathcal{F}_{T_n-}]\ind_{[T_n,\infty)}\\
&=&
\frac{1}{4^n}\mathbb{E}[\mathsf{a}\transp\mathsf{a}-\mathsf{a}\transp p-{p}\transp\mathsf{a}+{p}\transp{p}|\mathcal{F}_{T_n-}]\ind_{[T_n,\infty)}\\
&=&
\frac{1}{4^n}(\mathfrak{D}_{{p}}-{p}\transp{p})\ind_{[T_n,\infty)},
\dce
$$
where $\mathfrak{D}_{{p}}$ denotes the diagonal matrix of diagonal vector $p$. For any vector $a=(a_0,\ldots,a_k)$, if $$
0=\transp a (\mathfrak{D}_{{p}}-{p}\transp{p}) a 
= \mathbb{E}[(\transp a \mathsf{a}- \transp a{p})^2|\mathcal{F}_{T_n-}],
$$
necessarily $(\transp a \mathsf{a}- \transp a{p})^2=0$ or $a_h=\transp a{p}$ for all $1\leq h\leq k$. This means that the kernel of $(\mathfrak{D}_{{p}}-{p}\transp{p})$ is the vector space $\mathscr{K}$ generated by the vector $(1,1,\ldots,1)\in\mathbb{R}^k$, while its image space, as $(\mathfrak{D}_{{p}}-{p}\transp{p})$ is symmetric, is $\mathscr{K}^\perp\subset\mathbb{R}^k$. The matrix $(\mathfrak{D}_{{p}}-{p}\transp{p})$ as an operator on $\mathscr{K}^\perp$ is invertible. This implies the existence of an $\mathcal{F}_{T_n-}$ measurable matrix $\mathsf{J}$ such that $
\mathsf{J}(\mathfrak{D}_{{p}}-{p}\transp{p}),
$ 
on $\{\mathtt{I}_n=\{1,\ldots,k\}\}\cap\{T_n<\infty\}$, is the projection operator onto the space $\mathscr{K}^\perp$.

We can make the same analysis with $[\widetilde{X}'',\transp \widetilde{X}'']^{\mathbb{G}\cdot p}$. Notice that $\Delta_{T_n}\widetilde{X}''_h=\frac{1}{2^n}(\ind_{A_{n,h}}-\overline{p}_{n,h})\ind_{[T_n,\infty)}$. Notice that, on the set $\{\mathtt{I}_n=\{1,\ldots,k\}\}\cap\{T_n<\infty\}$, $\overline{p}_{n,h}=0$ for $h>k$. We obtain then that the vector $(1,1,\ldots,1)\in\mathbb{R}^k$ is in the kernel of the matrix$$
{\mathbf{M}}:=\left(\mathbb{E}[\Delta_{T_n}\widetilde{X}''_i\Delta_{T_n}\widetilde{X}''_j|\mathcal{G}_{T_n-}]\right)_{1\leq i,j\leq k},
$$
and by the symmetry, the image of the ${\mathbf{M}}$ is contained in $\mathscr{K}^\perp$. Now, for any vector $a\in \mathscr{K}$, $
{\mathbf{M}}a=0={\mathbf{M}} \mathsf{J}(\mathfrak{D}_{{p}}-{p}\transp{p})a,
$
while for any vector $a\in \mathscr{K}^\perp$, $
{\mathbf{M}}a={\mathbf{M}} \mathsf{J}(\mathfrak{D}_{{p}}-{p}\transp{p})a,
$
proving ${\mathbf{M}}={\mathbf{M}} \mathsf{J}(\mathfrak{D}_{{p}}-{p}\transp{p})$. Finally,
$$
\dcb
\Delta_{T_n}[\widetilde{X}'',\transp \widetilde{X}'']^{\mathbb{G}\cdot p}
&=&
\left(
\dcb
\mathbf{M},&0\\
\\
0,&0
\dce
\right)
=
\left(
\dcb
{\mathbf{M}} \mathsf{J}(\mathfrak{D}_{{p}}-{p}\transp{p}),&0\\
\\
0,&0
\dce
\right)\\

&=&
\left(
\dcb
\mathbf{M},&0\\
\\
0,&0
\dce
\right)
\left(
\dcb
\mathsf{J},&0\\
\\
0,&0
\dce
\right)
\left(
\dcb
\mathfrak{D}_{{p}}-{p}\transp{p},&0\\
\\
0,&0
\dce
\right)\\

&=&
\Delta_{T_n}[\widetilde{X}'',\transp \widetilde{X}'']^{\mathbb{G}\cdot p}
\left(
\dcb
\mathsf{J},&0\\
\\
0,&0
\dce
\right)4^n
\Delta_{T_n}[{X}'',\transp {X}'']^{\mathbb{F}\cdot p}
\dce
$$
on $\{\mathtt{I}_n=\{1,\ldots,k\}\}\cap\{T_n<\infty\}$.	On this set, define$$
J_n
=
\left(
\dcb
\mathsf{J},&0\\
\\
0,&0
\dce
\right)4^n.
$$
Now, making the above computation on the set $\{\mathtt{I}_n=\mathtt{B}\}\cap\{T_n<\infty\}$ for any no-empty subset $\mathtt{B}$ of $\{1,\ldots,d\}$ (instead of $\{1,\ldots,k\}$), we obtain an $\mathcal{F}_{T_n-}$ measurable matrix valued random variable everywhere defined $J_n$ such that 
$$
\Delta_{T_n}[\widetilde{X}'',\transp \widetilde{X}'']^{\mathbb{G}\cdot p}
=
\Delta_{T_n}[\widetilde{X}'',\transp \widetilde{X}'']^{\mathbb{G}\cdot p}
J_n
\Delta_{T_n}[{X}'',\transp {X}'']^{\mathbb{F}\cdot p}.
$$
The lemma is proved with$$
G=
\sum_{n=1}^{\mathsf{N}^a-}
\Delta_{T_n}[\widetilde{X}'',\transp \widetilde{X}'']^{\mathbb{G}\cdot p}J_n\
\ind_{[T_n]}. \ \ok
$$

\

\subsection{Drift of $X^\circ_h$}

For $1\leq h\leq \mathsf{n}'''$, Let $\mu$ denote the jump measure of $X^\circ_h$. Let $\nu$ (resp. $\overline{\nu}$) be the $\mathbb{F}$ (resp. $\mathbb{G}$) compensator of $\mu$. 

\bl\label{Hmn}
We have $\mathscr{G}^{}(\mathbb{F},\mu)\ind_{[0,T]}\subset\mathscr{G}^{}(\mathbb{G},\mu)$. Let $g\in\mathscr{G}(\mathbb{F},\mu)$. Then, on $[0,T]$, $\Gamma(g{_*}(\mu-{\nu}))$ is continuous and $$
g{_*}(\mu-\overline{\nu})
=
g{_*}(\mu-{\nu})-\Gamma(g{_*}(\mu-{\nu})).
$$
In particular, $x{_*}(\mu-\overline{\nu})=\widetilde{{X}}^\circ_h$ on $[0,T]$. 
\el 

\textbf{Proof.} According to Lemma \ref{totinacc}, the support set of $\mu$ avoids any $\mathbb{G}$ predictable stopping time $U$ on $[0,T]$ so that $\ind_{[U]}{_*}\overline{\nu}=0$ on $[0,T]$. This implies  $\mathscr{G}^{}(\mathbb{F},\mu)\ind_{[0,T]}\subset\mathscr{G}^{}(\mathbb{G},\mu)$ (cf. \cite[Definition 11.16]{HWY}). On the other hand, $X=g{_*}(\mu-{\nu})-\Gamma(g{_*}(\mu-{\nu}))$ is a $\mathbb{G}$ local martingale on $[0,T]$, whose jump $\Delta_UX$ at the $\mathbb{G}$ predictable stopping time $U$ is given by $-\Delta_U\Gamma(g{_*}(\mu-{\nu}))$ on $\{U\leq T\}$ which is $\mathcal{G}_{U-}$ measurable. Hence, $\Delta_UX=\Delta_U\Gamma(g{_*}(\mu-{\nu}))=0$ on $\{U\leq T\}$ (cf. \cite[Theorem 7.13]{HWY}), i.e. $\Gamma(g{_*}(\mu-{\nu}))$ is continuous on $[0,T]$. Now we compute the jumps on $[0,T]$.$$
\Delta_s(g{_*}(\mu-\overline{\nu}))
=
g(s,\Delta_s X^\circ_h)\ind_{\{\Delta_sX^\circ_h\neq 0\}}
=
\Delta_s(g{_*}(\mu-{\nu})-\Gamma(g{_*}(\mu-{\nu}))).
$$
The lemma is proved by \cite[Theorem 7.23]{HWY}. \ \ok

\bl\label{driftprime3}
There exists a $\mathbb{G}$ predictable process $G^\circ_{h}$ such that, on $[0,T]$,$$
\Gamma(X^\circ_h)
=
G^\circ_{h}\centerdot [X^\circ_h, X^\circ_h]^{\mathbb{F}\cdot p}.
$$
on $[0,T]$.
\el

\textbf{Proof.}
Let $Y$ be a $\mathbb{G}$ local martingale defined in Lemma \ref{fbd} for $X^\circ_h$. By Lemma \ref{projectionlemma}, there exists $g\in\mathscr{G}(\mathbb{G},\mu)$ such that, on $[0,T]$,
$$
[Y,X^\circ_h]^{\mathbb{G}\cdot p}
=
[Y,\widetilde{X}^\circ_h]^{\mathbb{G}\cdot p}
=
[g{_*}(\mu - \overline{\nu}), \widetilde{X}^\circ_h]^{\mathbb{G}\cdot p}.
$$
Using the notations of Section \ref{newmodif}, we verify that the time support of $\mu$ is $$
\{e_{h}(\alpha'''_{h})\ind_{\{\Delta X'''=\alpha'''_{h}\neq \boldsymbol{0}\}}\neq 0\}=\{\Delta X'''=\alpha'''_{h}\neq \boldsymbol{0}\},
$$
while its space location process can be $e_{h}(\alpha'''_{h})$. We see that $\mu$ satisfies the finite $\mathbb{F}$ predictable constraint condition (with constraint process $e_{h}(\alpha'''_{h})$) and satisfies the conditions in Theorem \ref{specif3} on $[0,T]$ (with $\mathsf{n}=1$). Hence, there exists a $\mathbb{G}$ predictable process $H$ such that
$$
[g{_*}(\mu - \overline{\nu}), \widetilde{X}^\circ_h]^{\mathbb{G}\cdot p}
=
[H\centerdot \widetilde{X}^\circ_h, \widetilde{X}^\circ_h]^{\mathbb{G}\cdot p}
=
H\centerdot [X^\circ_h, X^\circ_h]^{\mathbb{G}\cdot p}
$$
on $[0,T]$. By Lemma \ref{interm}, $[X^\circ_h, X^\circ_h]^{\mathbb{G}\cdot p}$ is absolutely continuous with respect to $[X^\circ_h, X^\circ_h]^{\mathbb{F}\cdot p}$ on $[0,T]$. Let$$
K=\frac{\mathsf{d}[X^\circ_h, X^\circ_h]^{\mathbb{G}\cdot p}}{\mathsf{d}[X^\circ_h, X^\circ_h]^{\mathbb{F}\cdot p}}.
$$
Then, on $[0,T]$,$$
\dcb
[Y,X^\circ_h]^{\mathbb{G}\cdot p}
&=&
H\centerdot [X^\circ_h, X^\circ_h]^{\mathbb{G}\cdot p}
=
HK\centerdot [X^\circ_h, X^\circ_h]^{\mathbb{F}\cdot p}.
\dce
$$
Together with Lemma \ref{fbd}, this concludes the proof.\ \ok

\

\

\end{document}